\documentclass[twocolumn]{autart}    

\usepackage{graphicx}          

\usepackage{amsmath,amssymb,mathtools,enumerate,enumitem}
\usepackage{subfigure}

\usepackage{algorithm,algpseudocode}

\usepackage{float}

\usepackage{natbib}

\newcommand\setb[1]{\left\{#1\right\}}
\newcommand{\real}{\mathbb{R}}

\newcommand{\abs}[1]{\left\vert #1 \right\vert}

\newcommand{\esc}[2]{\langle #1, #2 \rangle}
\newcommand{\norm}[1]{\left\lVert#1\right\rVert}
\DeclareMathOperator*{\argmin}{argmin}

\DeclareMathOperator{\prox}{Prox}
\DeclareMathOperator{\KL}{KL}

\graphicspath{ {figures/} }

\usepackage[dvipsnames]{xcolor}
\usepackage{tikz,pgfplots}
\usepackage{tkz-euclide}
\usepackage{pgfplotstable}
\usepackage{pgfplots}
\pgfplotsset{
	table/search path={plot_data},
	compat=1.17
}

\usepackage{xr}

\usepackage{environ}
\makeatletter
\newwrite\remember@figures
\AtBeginDocument{%
  \InputIfFileExists{\jobname.dft}{}{}%
  \immediate\openout\remember@figures=\jobname.dft
}
\AtEndDocument{\immediate\closeout\remember@figures}

\NewEnviron{dfigure}[1]{%
  \immediate\write\remember@figures{%
    \noexpand\rememberfigure{#1}{\unexpanded\expandafter{\BODY}}%
  }%
}
\newcommand{\placefigure}[2][tp]{%
    \csname remembered@figure@#2\endcsname{#1}
}
\newcommand{\rememberfigure}[2]{%
  \global\@namedef{remembered@figure@#1}##1{%
    \begin{figure}[##1]#2\end{figure}%
  }%
}

\NewEnviron{dfigurefull}[1]{%
  \immediate\write\remember@figures{%
    \noexpand\rememberfigurefull{#1}{\unexpanded\expandafter{\BODY}}%
  }%
}
\newcommand{\placefigurefull}[2][tp]{%
    \csname rememberedfull@figurefull@#2\endcsname{#1}
}
\newcommand{\rememberfigurefull}[2]{%
  \global\@namedef{rememberedfull@figurefull@#1}##1{%
    \begin{figure*}[##1]#2\end{figure*}%
  }%
}
\makeatother

\begin{document}

\begin{frontmatter}

\title{Approximate Wasserstein Attraction Flows for  \\ Dynamic Mass Transport over Networks\thanksref{footnoteinfo}} 

\thanks[footnoteinfo]{This paper was not presented at any IFAC 
meeting. Corresponding author: Ferran Arqu\'e. Tel. +34934015752. 
Fax +34934015750.}

\author[author1]{Ferran Arqu\'e}\ead{ferran.arque@upc.edu},             
\author[author2]{C\'esar A. Uribe}\ead{cauribe@rice.edu},                
\author[author1]{Carlos Ocampo-Martinez}\ead{carlos.ocampo@upc.edu}     

\address[author1]{Automatic Control Department, Universitat Politècnica de Catalunya, Institut de Robòtica i Informàtica Industrial
(CSIC-UPC), Llorens i Artigas, 4-6, 08028, Barcelona, Spain}  
\address[author2]{Department of Electrical and Computer Engineering, Rice University, Houston TX 77006}  


\begin{abstract}                          
This paper presents a Wasserstein attraction approach for solving dynamic mass transport problems over networks. In the transport problem over networks, we start with a distribution over the set of nodes that needs to be ``transported'' to a target distribution accounting for the network topology. We exploit the specific structure of the problem, characterized by the computation of implicit gradient steps, and formulate an approach based on discretized flows. As a result, our proposed algorithm relies on the iterative computation of constrained Wasserstein barycenters. We show how the proposed method finds approximate solutions to the network transport problem, taking into account the topology of the network, the capacity of the communication channels, and the capacity of the individual nodes. Finally, we show the performance of this approach applied to large-scale water transportation networks.
\end{abstract}

\end{frontmatter}

\section{Introduction}\label{sec:1intro}

The formalization of the optimal transportation problem by Gaspard Monge in 1781~\citep{monge} and Leonid Kantorovich~\citep{kant3} evolved into a whole branch of mathematics~\citep{gradientFlows,user,villani} called \textit{optimal transport theory}. Optimal transport (OT) is based on the computation of the distance between two objects, probability distributions in many cases, commonly referred to as the \textit{Wasserstein distance}~\citep{Vaserstein}. 
Optimal transport has a wide range of applications
, e.g., image retrieval~\citep{rubner}, 
averaging atmospheric gas concentration data sets~\citep{averagingGas}, 
segmentation and labeling of neurons~\citep{neuronLabel}, clustering patterns linking COVID-19 dynamics and human mobility~\citep{covid}, and different machine learning problems, such as generative adversarial networks~\citep{GAN}, low-rank matrix factorization~\citep{matrixFact} or fair regression~\citep{FairML}.

Despite the appealing properties of the Wasserstein distance, its computation requires solving an optimization problem. Such optimization problems become computationally prohibited on high dimensional objects, a large number of distributions or a high desired accuracy. However, the seminal work of Brenier~\citep{brenier} led to practical numerical algorithms that started the search for efficient algorithms to solve the OT problem~\citep{compOT}.

Our work focuses on the discrete OT problem, where probability distributions are defined over the nodes of a graph, assumed to be finite. 
In traditional OT approaches, it is assumed that mass (or a fraction of it) at each point in the support of one of the probability measures can be sent to any of the elements in the support of the other probability measure. As a result, the transport plan is executed effectively in one step. However, we seek to explicitly consider the topology of the underlying graph, which naturally imposes some transportation constraints. Adding the topology of the graph as a constraint means that there may not be a direct link between two points in the support, as the edges of the graph directly determine links. Therefore, our goal is to find a sequence of transport plans that move the mass from an initial distribution to a final one along the edges of a general connected graph so that the cost of transportation is minimal while accounting for channel and node capacities.

Finding the amount of mass that needs to be sent through each edge so that the total cost of transportation is minimal is a well-known problem called the \textit{minimum-cost flow problem}. This problem has been widely studied \citep{MCF1967,MCF1993,MCF2020}, and different algorithms have been proposed to solve it \citep{kovacs}. More importantly, the Wasserstein distance can be rewritten as a minimum-cost flow problem when considering a complete bipartite graph~\citep{bassetti2}. This can be extended to more general graphs if one considers the shortest path distance as the cost of sending a resource unit from one node to the other. In this case, if we compute the minimum-cost flow, then the mass sent from one node to the other is known. Moreover, since the shortest path between them is also known, we can make the first step of the transportation by sending that mass through the first edge of the path. Following this approach, we obtain the desired sequence of transport plans. However, classical methods to solve this problem do not have a condition to discern between paths when the optimal flow is not unique \citep{solomon2}. This nonuniqueness leads to unpredictability of the output from the solver since many paths can be indistinguishable in terms of costs. To avoid that case, some algorithms introduce an additional term to the objective function so that it becomes strongly convex. These regularized OT methods, like the well-known Sinkhorn algorithm \citep{cuturiSinkhorn}, achieve uniqueness and significantly speed up the computation, compared to solving a large linear programming problem. Still, it is at the cost of finding an approximation of the solution to the original problem.

Our approach is based on the resolution of the \textit{Wasserstein attraction} (WA) problem~\citep{PeyreGradientFlows}, which requires the computation of a \textit{Wasserstein barycenter} (WB) of two distributions at every iteration. Computing the WB yields an intermediate distribution, defined as the Fréchet mean of the two measures, which is the result of minimizing the sum of the (Wasserstein) distances between itself and each of the two distributions \citep{cuturiFastWB}. However, the support of this resulting distribution can include any of the graph nodes. We expand the definition of the WB problem by adding constraints that ensure the mean obtained has the appropriate support and each node does not receive more mass than the amount available from its neighbors. This approach resembles what is called \textit{displacement interpolation}~\citep{villani,solomonConvolutional}. However, displacement interpolation in the discrete-time case may require a small step size of the weight to ensure certain smoothness in the transportation (i.e., to avoid some of the mass moving over more than one node in a single step), which may lead to having many more iterations than necessary. Furthermore, with this approach, there is the possibility that certain nodes receive more mass than the total obtainable from their neighboring nodes. In summary, the main differentiating factor between displacement interpolation and our proposal is the addition of the topology and capacity constraints imposed by a graph. In this regard, \citep{Haasler2021} recently studied this problem in the context of traffic planning, where edge capacity constraints are taken into account, and proposed a framework based on the Lagrangian dual problem to solve it, which resembles the Sinkhorn algorithm.

Moreover, our proposed approach can be reformulated as a discrete gradient flow problem. Several papers work on discrete gradient flows over graphs (or other discrete domains) \citep{Chow2017ADS,erbar,Mielke2013GeodesicCO,Richemond2017OnWR}. 
However, such papers focus on the theoretical analysis of differential equations rather than the computational aspect with the regularized approximation of the Wasserstein metric (except for \citep{erbar} which provides a more in-depth discussion on the topic), and no additional constraints are considered on the elements of the graph. The closest works to our setting with constrained WB are \citep{PeyreGradientFlows,Cuturi2016ASD}. The former presents a framework to approximate gradient flows for Wasserstein metrics by computing discrete entropy-regularized flows, which are computed as \textit{JKO flows} (named after the authors in \citep{JKO}). It introduces the concept of Wasserstein attraction, which is used in our work. We expand on this concept by observing that our particular problem formulation allows us to write each iteration of the WA problem as the computation of a WB, which unlocks the use of powerful computational tools found in the literature to solve this problem. Additionally, as previously mentioned, we further generalize the definition of this regularized flow by including the supplemental constraints of the topology of a network and the node and edge capacity bounds, which are features not considered in~\citep{PeyreGradientFlows}.
The latter work, \citep{Cuturi2016ASD}, complements \citep{PeyreGradientFlows} while focusing on the dual formulation of Wasserstein variational problems. In the context of applications of JKO flows in OT, \citep{bunne2021} recently proposed a novel procedure for the computation of JKO flows based on input convex neural networks. It is applied in the study of population dynamics, where it assumes that the dynamics of the model is parameterized by an energy function, which controls how the transport is executed at each step, from one state to the next. In our application, this role is performed by another Wasserstein distance function instead of an energy one (in addition to further constraints), which also allows for explicit computation of the JKO steps. Similar works such as~\citep{Chen2018,Guex2019} propose methods for optimal transportation over networks based on Markov processes. The authors in~\citep{Chen2018} use the relative entropy as an index of closeness between measures and, doing so, they solve \textit{Schr\"odinger's bridge problem}~\citep{Schroedinger1931} for the computation of transport plans in a fixed number of steps. This same entropy is used to measure how much the mass spreads in the transportation. Thus, they design a transport plan where the mass spreads as much as possible to guarantee robustness against failures in the paths of the network, while still ensuring a reasonably low total cost. In a similar fashion, \citep{Guex2019} uses a bag-of-paths framework equivalent to solving either a standard or a relaxed entropy-regularized Wasserstein distance problem. Our approach allows topology changes as well, but it does so by solving a new problem at each step of the transportation, ignoring previous events, which increases the computational costs. 
However, the problem formulation allows us to consider node and edge capacity constraints explicitly. Moreover, flow speed can be adjusted with the weight parameter introduced when solving a WB problem.

The main contributions of this paper are threefold. First, we propose the mathematical formulation of a Wasserstein attraction-like problem to solve mass transport problems over networks by writing them as the computation of a WB problem with additional constraints. Second, we present a methodology to find an approximation of optimal discrete flows over networks based on Dykstra's projection algorithm and the computation of JKO flow proximal operators for the Kullback-Leibler divergence and prove the convergence of these intermediate steps under certain assumptions. Finally, to the best of our knowledge, there are no works related to water management systems under the Wasserstein distance framework. Hence, we illustrate how this approach using WB can be implemented to model a supply-and-demand problem in the context of drinking water networks, where the network constraints are a crucial aspect inherent in their nature. In addition, we show how it can automatically adapt to dynamic changes on the network's topology and agents. Furthermore, since there is no known method that can be used for fair comparison that can generate a flow that minimizes the Wasserstein distances and takes into account the network constraints, we have opted to compare the performance of our method with the commercial solver CPLEX with an explicit formulation of the constraints.

The remainder of this article is structured as follows. In Section \ref{sec:2prob_st}, we provide the necessary background for our work, stating some basic definitions from discrete OT theory and present the formal statement of the problem we want to solve. In Section \ref{sec:3proposed_approach}, we briefly review Dykstra's projection algorithm in the setting of optimization problems involving the Kullback-Leibler divergence and how it can be used to solve the WB problem. Then, we show the additional steps needed on the algorithm to enforce support constraints and capacity bounds on the network's links and nodes. With that, we present our proposed approach. In Section \ref{sec:4numerical_simulations}, we provide some illustrative examples. We discuss our approach in the context of flow optimization on drinking water networks and give some remarks regarding the numerical implementation of the proposed algorithm. Finally, in Section \ref{sec:5conclusions}, we provide some final comments and discuss future investigation directions.

\subsection*{Notation}

The column vector of all ones is denoted by $\boldsymbol{1}$ and $I$ is the identity matrix.
The adjacency matrix of a graph is denoted by $A$, and we will write $\bar{A} = A+I$ when considering the connection of one node to itself. $\real_+$ and $\real_{++}$ refer to non-negative and strictly positive real values respectively. Given $x\in\real^n$, $\norm{x}$ stands for its Euclidean norm. Given two matrices $A,B\in\real^{n\times m}$, $\esc{A}{B} = \sum_{i,j}A_{ij}B_{ij}$. We define the support of a function (or vector) $\rho$ as $\textsc{supp}(\rho) = \{i \mid \rho(i)>0 \}$.
We denote $\text{KL}(\pi|\xi)$ as the \textit{Kullback-Leibler divergence} between $\pi\in\real_+^{n\times n}$ and $\xi\in\real_{++}^{n\times n}$, defined as
\begin{equation*}
	\text{KL}(\pi|\xi) = \sum_{i,j=1}^n\pi_{ij}\ln\left(\frac{\pi_{ij}}{\xi_{ij}}\right) - \pi_{ij} + \xi_{ij},
\end{equation*}
with the convention $0 \ln(0) = 0$. Finally, the indicator function of a set $\mathcal{C}$ is defined as $\iota_\mathcal{C}(x) = 0$ if $x\in \mathcal{C}$, and $\iota_\mathcal{C}(x) = +\infty$ otherwise.
\section{Problem Statement: Discrete Flows and Wasserstein Attraction on Graphs}\label{sec:2prob_st}

\subsection{Discrete Flows on Graphs}\label{sec:2.1}

Consider a \textit{discrete}, \textit{finite}, \textit{fixed} and \textit{connected} graph $\mathcal{G} = (V,E)$, where $V$ is a set of $n$ nodes $V = (1,\cdots,n)$, and $E$ is a set of directed edges such that $E\subseteq V\times V$, where $(j,i) \in E$ if and only if there is an edge between the node $j \in V$ and node $i\in V$. Denote the probability simplex on $V$ as $\text{Prob}(V) = \{\mu \in \mathbb{R}_+^n  \mid \sum_{x\in V} \mu(x) =1\}$.
The set of edges $E$ has an associated weight function $c : E \to \mathbb{R}_+$ where each edge $e\in E$ has a corresponding weight $c_e = c(e)$, i.e., the cost of sending a unit of mass using the edge $e$. Furthermore, endow the graph $\mathcal{G}$ with its natural metric $\mathsf{d}$ which measures the total weight of the shortest path between any two nodes in $\mathcal{G}$. 

We study the discrete flow (i.e., discretization in time) problem of optimally transporting an initial mass distribution $\mu \in \text{Prob}(V)$ to a target mass distribution $\nu \in \text{Prob}(V)$ using the graph $\mathcal{G}$. The associated weight of each edge allows us to define a cost matrix $C \in \mathbb{R}_{+}^{n\times n}$, where $[C]_{ji} = \mathsf{d}(j,i)$ indicates the cost of transporting a unit mass from node $j$ to node $i$. Moreover, we endow the space $\text{Prob}(V)$ of probability measures on $V$ with the \mbox{$1$-Wasserstein} distance between two probability distributions $\mu$ and $\nu$ on $\mathcal{G}$ as
\begin{align*}
	W_1(\mu,\nu) & = \min_{\pi \in \Pi(\mu,\nu)} \sum_{x,y \in V} \mathsf{d}(x,y) \pi(x,y),
\end{align*}
where the minimizer is computed over all couplings on $V\times V$ with marginals $\mu$ and $\nu$, i.e., the set of \textit{optimal transport plans} $\Pi(\mu,\nu) = \setb{\pi\in\real^{n\times n}_+ \big| \pi\boldsymbol{1} {=} \mu, \pi^\intercal\boldsymbol{1}{=}\nu}$.

Our objective is to design a discrete flow $\{\rho_t\}_{t\geq 0}$ on $\mathcal{G}$, where $\rho_t \in \text{Prob}(V)$, by constructing a sequence of transport plans $\{\pi_t\}_{t\geq0}$ such that $\rho_0 = \mu$, $\rho_{t+1} = \pi_t \boldsymbol{1}$, $\rho_{t} = \pi_t^\intercal \boldsymbol{1}$ and $\lim_{t\to \infty} \rho_t = \nu$. Moreover, the transport cost at each iteration should be minimized.

Furthermore, the desired sequence of transport plans is required to satisfy the following constraints imposed by the network:
\begin{enumerate}[label=(\alph*)]
	\item A node can only send mass to its neighbors, i.e.,
	\begin{align} \label{cons:a}
		[\pi_{{t}}]_{ij} > 0 \ \ \text{if} \ \ [\rho_t]_j >0  \ \text{and} \ (j,i)\in E. 
	\end{align}
	In other words, the flow should follow the sparsity pattern induced by the graph topology. Intuitively, a flow can only be assigned between two nodes if and only if there is an edge connecting them. Hence, for a transport plan $\pi_t$ it must hold that \mbox{$\textsc{supp}(\rho_{t+1}) \subseteq \{\textsc{supp}(\rho_t) \cup \{ j \mid (j,i)\in E \} \}$}.
	\item The mass sent over an edge cannot be greater than the associated edge capacity, i.e.,
	\begin{align} \label{cons:b}
		\pi_t \leq \tilde{C},
	\end{align}
	for a matrix of capacities $\tilde{C} \in \mathbb{R}_{+}^{n\times n}$, where $[\tilde{C}]_{ij}$ is the capacity of the edge $(j,i) \in E$ (the inequality should be understood element-wise).
	\item The mass at a node $i$ at some time instant $t \geq 0$ must not exceed its local storage capacity, i.e.,
	\begin{align} \label{cons:c}
		\rho_t \leq \rho ,
	\end{align}
	for a vector of storage capacities $\rho \in \mathbb{R}_{+}^n$ (again, the inequality is understood entry-wise). 
	\item The mass transported from a node $j$ to a node $i$ cannot exceed the mass held at node $j$, i.e.,
	\begin{align*}
		[\pi_{{t}}]_{ij} \leq [\rho_t]_j.
	\end{align*}
\end{enumerate}

\subsection{Wasserstein Attraction Flows}

We formulate the dynamic transport problem described in Section \ref{sec:2.1} as a constrained \textit{Wasserstein attraction} (WA) problem~\citep[Section 5.2]{PeyreGradientFlows}. Our main technical tool will be the JKO flow proximal operators which we introduce next.  We first present the JKO flow proximal operator with respect to a functional $f$. For all $q \in \text{Prob}(V)$,
\begin{align*}
	\text{Prox}^{W_1}_{\tau,f}(q) \triangleq \argmin_{p \in \text{Prob}(V)} \left\lbrace W_1(p,q) + \tau f(p) \right\rbrace,
\end{align*}

where $\tau$ is a step-size. Thus, starting from an initial distribution $\rho_0 = \mu$, the discrete JKO flow with respect to $f$ is defined as
\begin{align}\label{eq:implGradSteps}
	\rho_{t+1} \triangleq \text{Prox}^{W_1}_{\tau,f}(\rho_{t}).
\end{align}

Wasserstein attraction refers to the flow generated by the implicit gradient steps in \eqref{eq:implGradSteps}, known as \textit{JKO stepping}, with respect to the potential function defined as $W_1(\rho_t,\nu)$ for some fixed distribution $\nu$. Informally, the potential function drives the flow to minimize its Wasserstein distance to a \textit{target} distribution. Thus, we define the WA discrete flow as
\begin{align}\label{eq:nonreg_bar}
	\rho_{t+1} & = \text{Prox}^{W_1}_{\tau,W_1(\cdot,\nu)}(\rho_{t}) \nonumber \\
	&  = \argmin_{p \in \text{Prob}(V)} \left\lbrace W_1(p,\rho_t) + \tau W_1(p,\nu) \right\rbrace.
\end{align}

The WA defined in~\eqref{eq:nonreg_bar} has a precise optimization structure. However, the computation of each proximal operation is computationally intense~\citep{PeyreGradientFlows}. Moreover, the constraints imposed by the graph are not taken into account. In the next subsection, we describe our proposed approach for the efficient computation of the discrete WA, taking into account the constraints imposed by the network.

\subsection{Approximate Wasserstein Attraction Flow on Graphs}

Initially, we present the entropy regularized discrete JKO flow for the WA problem following the ideas introduced in~\citep{PeyreGradientFlows}. The main contribution in~\citep{PeyreGradientFlows} is to replace the Wasserstein distance functions with their entropy regularized versions. The use of entropic regularization has been shown useful for the design of computational approaches for OT \citep{cuturiSinkhorn}.

\begin{defn}
	Given a cost matrix $C\in\real^{n\times n}_+$, the \textbf{discrete entropy-regularized Wasserstein distance} between $\mu,\nu \in \text{Prob}(V)$ is defined as
	\begin{equation}\label{eq:wass_dist}
		\mathcal{W}_\gamma(\mu,\nu) = \min_{\pi\in\Pi(\mu,\nu)} \esc{C}{\pi} + \gamma H(\pi),
	\end{equation}
	where $H(\pi) = \sum\pi_{ij}(\ln\pi_{ij}-1) = \esc{\pi}{\ln\pi - \boldsymbol{11^\top}}$ is the negative entropy and $\gamma\geq0$ is the regularization parameter.
\end{defn}

Now, we can define the approximate entropy-regularized WA flow as
\begin{align}\label{eq:reg_bar}
	\rho_{t+1} & = \text{Prox}^{\mathcal{W}_\gamma}_{\tau,\mathcal{W}_\gamma(\cdot,\nu)}(\rho_{t}) \nonumber \\
	&  = \argmin_{p \in \text{Prob}(V)} \left\lbrace \mathcal{W}_\gamma(p,\rho_t) + \tau \mathcal{W}_\gamma(p,\nu) \right\rbrace.
\end{align}
Note $\mathcal{W}_\gamma(\cdot,\cdot)$ is a strictly convex and coercive function, therefore the operator in~\eqref{eq:reg_bar} is uniquely defined.

Next, we state one important observation about the entropy-regularized WA flow in~\eqref{eq:reg_bar}. Without loss of generality, one can multiply the argument in the optimization problem~\eqref{eq:reg_bar} by a constant $\omega = 1/(1+\tau)$. Thus, we obtain
\begin{align}\label{eq:reg_bar2}
	\rho_{t+1} & = \argmin_{p \in \text{Prob}(V)} \left\lbrace \omega \mathcal{W}_\gamma(p,\rho_t) + (1-\omega) \mathcal{W}_\gamma(p,\nu) \right\rbrace,
\end{align}
which is precisely the entropy-regularized Wasserstein barycenter between $\rho_{t}$ and $\nu$ \citep{cuturiFastWB}. Recall that for a finite set of probability distributions $\{\mu_i\}_{i=1}^m$ where \mbox{$\mu_i \in \text{Prob}(V)$}, the entropy-regularized Wasserstein barycenter is defined as
\begin{align*}
	\mu \triangleq \argmin_{p \in \text{Prob}(V)} \sum_{i=1}^m \omega_i \mathcal{W}_\gamma(p,\mu_i),
\end{align*}
where $\omega_i \geq 0$ and $\sum_{i=1}^m \omega_i =1$.

We interpret the Wasserstein attraction problem as the sequential computation of Wasserstein barycenters. This introduces an additional weight parameter that can be modified to give preference to one measure or the other. Such parameter consequently alters how the mass is transported across the graph, as we illustrate further along this paper.

Note that the barycenter is not restricted to only two distributions but as many as one may need. This means that the solution proposed here could be extended for problems akin to ours but involving more than two distributions, and in turn, we would have several weight parameters to customize the solutions obtained~\citep{multimarginalOT}.

The method that we propose uses \textit{Dykstra’s projection algorithm} \citep{Dykstra}. In our setting, much like Sinkhorn's algorithm, it is easier to implement than more traditional schemes designed to solve mathematical programs.

Another feature of the proposed approach is that, unlike in the computation of the Wasserstein distance (or, for that matter, solving the minimum-cost flow problem), we do not compute the complete flow in a single step, which would also entail having to store the shortest path between each node (or at least the first step of each path). In this regard, our method not only does not need to store this additional information, but it is also \textit{memoryless} in the sense that, at each step, the algorithm solves a new problem with initial and final distributions. This is advantageous since these measures do not need to be the same as in the previous steps (even the parameters, such as the weights, can be changed). This adaptability is the main difference between the flow we compute, a discrete one, and the one found by solving a minimum-cost flow problem, which is continuous. These aspects might take importance in future works where this method could be adapted in the context of decentralized or distributed optimization, where the available information at each node is limited~\citep{distributedQuantization,decentralizeAndRand}.

Approximate solutions to problems of the form~\eqref{eq:reg_bar2} can be efficiently computed by reformulating the entropy-regularized OT problem~\eqref{eq:wass_dist} as
\begin{equation}\label{eq:IBP}
	W_\gamma(\mu,\nu) = \min_{\pi\in \Pi(\mu,\nu)} \text{KL}(\pi|\xi),
\end{equation}
where $\xi = e^{-C/\gamma}$ (entry-wise exponential)~\citep{IBP}. Note that~\eqref{eq:IBP} can be extended for higher dimensional arrays (such as the tuples $\boldsymbol{\pi} = (\pi_1,\ldots,\pi_m)$ introduced in the definition of the WB) by summing over the indices $(i,j,k,\ldots)$. Thus, following~\citep{IBP}, we can rewrite Problem~\eqref{eq:reg_bar2} as
\begin{align}\label{eq:wasb_kl}
	\min_{\boldsymbol{\pi}\in \mathcal{C}_f\cap\mathcal{C}_e} \text{KL}_\omega(\boldsymbol{\pi}| \xi) {=} \omega \text{KL}(\pi_1|\xi) {+} (1{-}\omega)\text{KL}(\pi_2|\xi),
\end{align}
where
\begin{align}
    \mathcal{C}_f &= \setb{ \pi_1,\pi_2 \mid \pi_1\boldsymbol{1} = \rho_t, \pi_2\boldsymbol{1} = \nu},\label{eq:C1}\\[7pt]
    \mathcal{C}_e &= \setb{\pi_1,\pi_2 \mid \pi_1^\intercal \boldsymbol{1} =   \pi_2^\intercal \boldsymbol{1} = p}.\label{eq:C2}
\end{align}

Finally, taking into account the constraints in~\eqref{cons:a},~\eqref{cons:b} and~\eqref{cons:c} in Problem~\eqref{eq:wasb_kl}, we can state our main contribution regarding the design of the entropy-regularized discrete WA flow.

\begin{prob}\label{prob:prob_statement}
Consider a discrete, finite, fixed and connected graph with $n$ vertices, $\Tilde{C}\in\real_+^{n\times n}$ the capacity matrix, and $\mu,\,\nu\in \text{Prob}(V)$ the initial and final distributions respectively. We design the sequence of probability measures $\{\rho_t\}_{t\geq 0}$ by finding, for each $t \geq 0$, the transport plan that solves the optimization problem 
	\begin{subequations}\label{eq:prob_st}
		\begin{align}
			&\{\pi_t\} {=} \argmin_{\substack{\boldsymbol{\pi}\in\mathcal{C}_f \cap \mathcal{C}_e\\ \boldsymbol{\pi}\in\mathcal{C}_1 \cap\mathcal{C}_2 \cap \mathcal{C}_3 }} \omega \text{KL}(\pi_1|\xi) {+} (1{-}\omega)\text{KL}(\pi_2|\xi), \label{eq:prob_st_1}\\
			\intertext{where}
			&\mathcal{C}_f = \setb{\boldsymbol{\pi}\in\real^{n\times n}_{+}{\times}\real^{n\times n}_{+}  \mid \pi_1\boldsymbol{1} = \rho_t,\pi_2\boldsymbol{1} = \nu}  \label{eq:prob_st_2}\\[3pt]
			&\mathcal{C}_e = \setb{\boldsymbol{\pi}\in\real^{n\times n}_{+}{\times}\real^{n\times n}_{+}  \mid \pi_1^\intercal \boldsymbol{1} = \pi_2^\intercal\boldsymbol{1} = p}  \label{eq:prob_st_3}\\
			&\mathcal{C}_1 = \setb{\boldsymbol{\pi}\in\real^{n\times n}_{+}{\times}\real^{n\times n}_{+} \mid \pi_1 \leq  \tilde C } \label{eq:prob_st_4}\\
			&\mathcal{C}_2 = \setb{\boldsymbol{\pi}\in\real^{n\times n}_{+}{\times}\real^{n\times n}_{+} \mid \pi_1^\intercal \boldsymbol{1} \leq  \rho,\pi_2^\intercal \boldsymbol{1} \leq  \rho}   \label{eq:prob_st_5}\\
			&\mathcal{C}_3 = \bigg\{\boldsymbol{\pi}\in\real^{n\times n}_{+}{\times}\real^{n\times n}_{+} {\mid} [\pi_1^\intercal \boldsymbol{1}]_i {\leq}  {\sum_{j : (j,i) \in E}} [\rho_t]_j \bigg\}  \label{eq:prob_st_6}
		\end{align}
	\end{subequations}
\end{prob}

We note that the constraint set $\mathcal{C}_3$ is redundant if in $\mathcal{C}_1$ we consider $[\Tilde{C}]_{ij}=0$ when nodes $i$ and $j$ are not connected. This is, in fact, what we propose for our procedure in Section~\ref{sec:descr_proposed_aproach}. Nevertheless, we write it explicitly in the problem formulation since it is a necessary constraint that could be imposed differently in other methodologies.

Figure \ref{fig:example_prob_st} shows a simple example to illustrate the steps we obtain by solving Problem \ref{prob:prob_statement}. The transport is computed over a path graph, and it starts with an initial distribution (top left) with its mass located in the central nodes, and the final mass (bottom right) is distributed closer to the extremes of the path. At each iteration, we show the resulting distribution found by solving \eqref{eq:prob_st} considering the previous solution as the initial measure. We have also considered a storage capacity of $0.3$ for the third-to-last node, resulting in partially sending the mass in the fourth iteration. We see how the mass is transported from the initial setting until the final distribution is eventually covered while verifying the constraints imposed in the problem statement addressed here.

\placefigurefull[t]{Figure1}

\section{Iterative Projections for the Computation of Transport Plans}\label{sec:3proposed_approach}

Now that we have the necessary background on discrete OT and have introduced the problem we want to solve, we describe the approach that we propose. We will solve the regularized version of the WB problem, with the additional constraints \eqref{cons:a}, \eqref{cons:b} and \eqref{cons:c}. To do so, we use a well-known algorithm for solving regularized OT problems called Dykstra's projection algorithm \citep{Dykstra}, which, in our setting, is a generalization of the widely used Iterative Bregman Projections (IBP) algorithm \citep{IBP}. We use Dykstra's method because the convergence of IBP cannot be guaranteed in the presence of inequality constraints.

In Section~\ref{sec:IBP_Dykstra}, we give some background on how this algorithm is used to compute the regularized WB. In Section~\ref{sec:cap_supp_constr}, we show how one can modify the algorithm to compute the WB with the added constraints, and finally, in Section \ref{sec:descr_proposed_aproach}, we move on to the description of the proposed algorithm.

\subsection{Computation of the WB using Dykstra's projection algorithm} \label{sec:IBP_Dykstra}

Dykstra's projection algorithm can be used to solve problems of the form
\begin{equation*}
	\min_{\pi\in\cap_i\mathcal{C}_i} \text{KL}(\pi|\xi),
\end{equation*}
much like Problem \ref{prob:prob_statement} defined in Section \ref{sec:2prob_st}. It is based on the computation of the proximal operators for the KL divergence. This is done iteratively, cycling through each constraint set $\mathcal{C}_i$, and since $\mathcal{C}= \cap_i\mathcal{C}_i$ is a finite intersection of $L$ sets, we shall define, for every index $i$, $\mathcal{C}_{i+L} = \mathcal{C}_i$. Then, for each $k>0$ we compute
\begin{equation*}
	\pi^{(k)} = \prox_{\iota_{\mathcal{C}_k}}^{\KL}\left(\pi^{(k-1)}\cdot q^{(k-L)}\right), \,\, q^{(k)} = q^{(k-L)}\frac{\pi^{(k-1)}}{\pi^{(k)}},
\end{equation*}
with initial values $\pi^{(0)} = \xi \quad\text{ and }\quad q^{(0)} = q^{(-1)} = \ldots = q^{(-L+1)} = \boldsymbol{1}\boldsymbol{1}^\intercal$. The product and division of matrices are considered element-wise. We slightly abuse notation by omitting the step-size $\tau$ in the definition of the proximal operator, since we are multiplying the argument in the optimization problem \eqref{eq:reg_bar} by $\omega = 1/(1+\tau)$, as noted in Section \ref{sec:2prob_st}.

The next proposition states how we can compute in closed form the proximal operator corresponding to each constraint in the WB problem~\eqref{eq:wasb_kl}.

\begin{prop}[Proposition 1 in \citep{IBP}]\label{prop:prox_C1}
	The proximal operator of the indicator function $\iota_{\mathcal{C}_f}$, corresponding to the constraint set $\mathcal{C}_f$ in \eqref{eq:C1}, has the closed form
    \begin{equation}\label{eq:ProxC1}
		{\left[\prox_{\iota_{\mathcal{C}_f}}^{\KL_\omega}(\boldsymbol{\pi})\right]_l} {=} \prox_{\iota_{\{\pi_l\boldsymbol{1}{=}P_l\}}}^{\KL}\hspace*{-0.1cm}(\pi_l) {=} \emph{diag}{\left(\frac{P_l}{\pi_l\boldsymbol{1}}\right)}{\pi_l},
	\end{equation}
where $l=1,2$ and $P_1 = \rho_t,P_2=\nu$.
\end{prop}

\begin{rem}
	For set $\mathcal{C}_f$, since the constraint is imposed to each transport plan independently from the rest, we can compute the proximal operator $\prox_{\iota_{\mathcal{C}_f}}^{\KL_\omega}(\boldsymbol{\pi})$ the same way as with the Wasserstein distance in \citep{IBP}, but individually for each $\pi_l$.
\end{rem}

\begin{prop}[Proposition 2 in \citep{IBP}]\label{prop:prox_C2}
	The proximal operator of the indicator function $\iota_{\mathcal{C}_e}$, corresponding to the constraint set $\mathcal{C}_e$ in \eqref{eq:C2}, has the closed form
	\begin{equation}\label{eq:ProxC2}
		\left[\prox_{\iota_{\mathcal{C}_e}}^{\KL_\omega}(\boldsymbol{\pi})\right]_l = \pi_l\emph{diag}\left(\frac{p}{\boldsymbol{1}^\intercal\pi_l}\right),
	\end{equation}
	where ${p = \prod_{l=1}^m\left(\boldsymbol{1}^\intercal\pi_l\right)^{\omega_l}}$ (the products and exponentiation are considered element-wise), and $m=2$ in our case. 
\end{prop}

\subsection{Capacity and support constrained WB} \label{sec:cap_supp_constr}

In the context of networks, it is reasonable to restrict how much mass can be sent from one node to another, i.e., to add a capacity to the edges connecting the nodes. This constraint is imposed on each transport plan by defining a capacity matrix $\Tilde{C}{\in}\real^{n\times n}$ such that $[\Tilde{C}]_{ij}$ is the maximum mass that can be sent from node $i$ to node~$j$.

Similarly to Proposition \ref{prop:prox_C1}, since this capacity constraint is imposed on each transport plan independently of the rest, the projection is done 
individually for each transport plan. The following proposition concerns the computation of the proximal operator for the set $\mathcal{C}_1$ in~\eqref{eq:prob_st_4}.

\begin{prop}[Section 5.2 in \citep{IBP}]\label{prop:prox_C3}
	The proximal map for the function $\iota_{\{\pi_1\leq \Tilde{C}\}}$ is defined as
	\begin{equation}\label{eq:proxC3}
		\prox_{\iota_{\{\pi_1\leq \Tilde{C}\}}}^{\KL}(\pi_1) =  \min{\left(\pi_1,\Tilde{C}\right)},
	\end{equation}
	with the minimum computed element-wise.
\end{prop}

We can also have capacity limits on some of the nodes, meaning that even though the optimal solution might send a certain amount of mass to one of these nodes, it may not be possible to hold that much quantity. 
This corresponds to the constraint set $\mathcal{C}_2$ in~\eqref{eq:prob_st_5}.
This set is, in fact, the same as one of the sets defined to solve partial transport problems, as seen in \citep{IBP} (except for having $\rho$ instead of one of the marginals). From that, we get the following result for the computation of the projection on this set in closed form.

\begin{prop}[Proposition 5 in \citep{IBP}]\label{prop:prox_C4}
	For the the indicator function $\iota_{\mathcal{C}_2}$, corresponding to the constraint set $\mathcal{C}_2$ in \eqref{eq:prob_st_5}, one has
	\begin{equation}\label{eq:proxC4}
		\begin{aligned}
			\left[\prox_{\iota_{\mathcal{C}_2}}^{\KL_\omega}(\boldsymbol{\pi})\right]_l &= \prox_{\iota_{\{\pi_l^\intercal\boldsymbol{1}\leq \rho\}}}^{\KL}(\pi_l)\\[5pt]
			&= \pi_l\emph{diag}\left(\min\left(\frac{\rho}{\pi_l^\intercal\boldsymbol{1}},\boldsymbol{1}\right)\right),
		\end{aligned}
	\end{equation}
	where the minimum and division of vectors are considered element-wise.
\end{prop}

In addition to the capacity constraints \eqref{eq:prob_st_4} and \eqref{eq:prob_st_5}, we want to restrict the barycenter domain to a smaller set of nodes, rather than the whole graph, since we can only send mass to the nodes in the support of $\rho_t$ and their neighbors. In this case, we would obtain a vector of dimension $n^*\leq n$, where each element corresponds to the mass at one node of the subset. It is clear from the second constraint of the WB problem, $\pi_l^\intercal\boldsymbol{1} = p \,\,\forall l$ for some measure $p$, that by resizing $p$ to have dimension $n^*$, now $\pi_l\in\real^{n\times n^*}$ and, subsequently, for the cost and capacity matrices, we should only take the columns corresponding to the subset of nodes (thus $C_l,\,\Tilde{C}_l\in\real^{n\times n^*}$). Therefore, the dimensions of the arrays in the computation of the projections still agree. However, let us go through the deduction of the computation of the projection on $\mathcal{C}_e$ shown in \eqref{eq:ProxC2} to see that it is well defined and it still holds for this support constraint ($\prox_{\iota_{\mathcal{C}_f}}^{\KL_\omega}$ is similar and $\prox_{\iota_{\mathcal{C}_1}}^{\KL_\omega}$ and $\prox_{\iota_{\mathcal{C}_2}}^{\KL_\omega}$ are straightforward).

\begin{prop}
	The computation of the proximal operator $\prox_{\iota_{\mathcal{C}_e}}^{\KL_\omega}(\boldsymbol{\pi})$ in \eqref{eq:ProxC2} still holds for $n\times n^*$ dimensional matrix inputs, where $n^*\leq n$.
\end{prop}
\begin{pf}
	Given $\pi_l^{(k-1)}\in\real^{n\times n^*}$, computing the projection on the set $\mathcal{C}_e$ consists in solving the optimization problem
	\begin{equation*}
		\min_{\boldsymbol{\pi}^{(k)}\in\mathcal{C}_e} \sum_{l=1}^m \omega_l\text{KL}\left(\pi_l^{(k)}\Big|\pi_l^{(k-1)}\right).
	\end{equation*}
	
	For the sake of notation, we define $\pi_l \coloneqq \pi_l^{(k)}$, $\overline{\pi}_l \coloneqq \pi_l^{(k-1)}$ and, with that, expanding the problem leaves us with
	\begin{align*}
		&\min_{\boldsymbol{\pi},p} \sum_{l,i,j} \omega_l\pi_{l,ij}\left(\ln\frac{\pi_{l,ij}}{\overline{\pi}_{l,ij}} - 1\right)\\
		&\text{s.t. } \, \pi_l^\intercal \boldsymbol{1} = p, \qquad l=1,\ldots,m.
	\end{align*}
	
	The Lagrangian of this problem is
	\begin{equation*}
		\sum_l\sum_{i,j} \omega_l \pi_{l,ij}\left(\ln\frac{\pi_{l,ij}}{\overline{\pi}_{l,ij} - 1}\right) + \lambda_l^\intercal\left(\pi_l^\intercal\boldsymbol{1}-p\right),
	\end{equation*}
	where $\lambda_l\in\real^{n^*}\, \forall l$ are the Lagrange multipliers.
	
	On one hand, if we differentiate the Lagrangian with respect to $\pi_{l,ij}$ and equate to zero, we get
	\begin{equation}\label{eq:PC2_lagr_cond1}
		\omega_l\left(\ln\frac{\pi_{l,ij}}{\overline{\pi}_{l,ij}}\right) + \lambda_{l,j} = 0.
	\end{equation}
	
	On the other hand, differentiating with respect to $p_j$ yields
	\begin{equation}\label{eq:PC2_lagr_cond2}
		-\sum_l \lambda_{l,j} = 0.
	\end{equation}
	
	Isolating $\pi_{l,ij}$ from \eqref{eq:PC2_lagr_cond1} yields $\pi_{l,ij} = \overline{\pi}_{l,ij}e^{\frac{-\lambda_{l,j}}{\omega_l}}$, thus
    \begin{equation}\label{eq:PC2_rel1}
		\pi_l = \overline{\pi}_l\text{diag}\left(e^{\frac{-\lambda_{l,1}}{\omega_l}}, \ldots, e^{\frac{-\lambda_{l,n^*}}{\omega_l}}\right).
	\end{equation}
	
	Then, combining \eqref{eq:PC2_rel1} together with the constraint $\pi_l^\intercal\boldsymbol{1} = p$, we obtain
	\begin{equation}\label{eq:PC2_rel2}
		\text{diag}\left(e^{\frac{-\lambda_{l,1}}{\omega_l}},\ldots, e^{\frac{-\lambda_{l,n^*}}{\omega_l}}\right)\overline{\pi}_l^\intercal\boldsymbol{1} = p,
	\end{equation}
	from which we deduce
	\begin{equation}\label{eq:PC2_rel3}
		\text{diag}\left(e^{\frac{-\lambda_{l,1}}{\omega_l}}, \ldots, e^{\frac{-\lambda_{l,n^*}}{\omega_l}}\right) = \text{diag}\left(\frac{p}{\overline{\pi}_l^\intercal\boldsymbol{1}}\right),
	\end{equation}
	where the division is considered element-wise.
	
	We still have to use the result in \eqref{eq:PC2_lagr_cond2}, so, we first rewrite \eqref{eq:PC2_rel2} as $\left(\overline{\pi}_l^\intercal\boldsymbol{1}\right)^{\omega_l} = \text{diag}\left(e^{\lambda_{l,1}}, \ldots, e^{\lambda_{l,n^*}}\right) p^{\omega_l}$,
	with element-wise exponentiation. With this relation, we compute $\prod_l \left(\overline{\pi}_l^\intercal\boldsymbol{1}\right)^{\omega_l}$, which is
	\begin{equation*}
		\prod_l \left(\overline{\pi}_l^\intercal\boldsymbol{1}\right)^{\omega_l} = \text{diag}\left(e^{\text{\normalsize$\Sigma$}_l\lambda_{l,1}}, \ldots, e^{\text{\normalsize$\Sigma$}_l\lambda_{l,n^*}}\right) p^{\text{\normalsize$\Sigma$}_l\omega_l}.
	\end{equation*}
	
	Then, using \eqref{eq:PC2_lagr_cond2} and the fact that $\sum_l\omega_l = 1$, we can write the measure $p$ in terms of the known quantities $\overline{\pi}_l$ and $\omega_l$ as
	\begin{equation}\label{eq:PC2_rel4}
		p = \prod_l \left(\overline{\pi}_l^\intercal\boldsymbol{1}\right)^{\omega_l}.
	\end{equation}
	
	Finally, combining \eqref{eq:PC2_rel1}$-$\eqref{eq:PC2_rel3}$-$\eqref{eq:PC2_rel4}, we obtain
	\begin{equation*}
		\pi_l = \overline{\pi}_l\text{diag}\left(\frac{p}{\overline{\pi}_l^\intercal\boldsymbol{1}}\right), \qquad \text{where } p = \prod_l \left(\overline{\pi}_l^\intercal\boldsymbol{1}\right)^{\omega_l},
	\end{equation*}
	which is what we wanted to show. 
	\qed
\end{pf}

\subsection{Description of the proposed approach} \label{sec:descr_proposed_aproach}

Now, we can present the proposed algorithm to solve Problem \ref{prob:prob_statement}. We use Dykstra's projection algorithm, and together with the support and capacity constraints, we can impose the additional restrictions introduced in the problem statement (Section \ref{sec:2prob_st}).

For the support constraint \eqref{eq:prob_st_6}, we will take for each matrix only the columns corresponding to the nodes in the support of $\rho_t$ and their neighbors, which we know, since we have the adjacency matrix $A$. Once we compute $\rho_{t+1}$, as it might have a smaller dimension $n^*\leq n$, we can redefine $\rho_{t+1}$ as an $n$-dimensional vector of all zeros except for the nodes that the elements of $\rho_{t+1}$ referred to, which will have the value that we have just computed. While this definition reduces the result to the desired support, nodes in $\textsc{supp}(\rho_t)$ can still send mass to non-neighboring nodes. To fix this issue, we adapt constraint $\eqref{eq:prob_st_4}$. We redefine the capacity matrix $\Tilde{C}$ for the transport plan $\pi_1$ from $\rho_t$ to $\rho_{t+1}$, such that for the nodes in the support of $\rho_t$, if there is no connection between one of them and another node, the "link" between them has zero capacity, i.e.,
\begin{equation}\label{eq:cap_mat}
	[\Tilde{C}]_{ij} = \begin{cases} 0 \quad &\text{if } j \in \textsc{supp}(\rho_t) \text{ and } \Bar{A}_{ij} = 0, \\ [\Tilde{C}]_{ij} &\text{otherwise}.\end{cases}
\end{equation}
We note that, in this case, constraints \eqref{eq:prob_st_4} and \eqref{eq:prob_st_6} could have had a separate matrix for each one and be considered two different projections on the algorithm, but here we merge both into one.

Resizing the matrices to limit the support is unnecessary for the algorithm to converge to the desired solution since it is already taken care of by the capacity matrix \eqref{eq:cap_mat}. Nevertheless, by implementing it, the dimension of the problem can be reduced, so the computations can be executed faster. In the worst-case scenario where all the nodes have mass, there is a direct connection to all the nodes or similar settings, the matrices and vectors would not be modified, and the algorithm would proceed as if this support constraint was not implemented.

Finally, for the storage capacity on each node \eqref{eq:prob_st_5}, we resize $\rho$ to only consider the elements corresponding to the nodes on the new support.

Algorithm \ref{alg:1} summarizes the proposed method. It is important to remark that our entropy-regularized approach does not allow the scheme to converge exactly to the target distribution $\nu$. Since the additional entropy term in the definition of the Wasserstein distance \eqref{eq:wass_dist} forces every node to send a small amount of mass to the rest, even if it does not correspond to the distribution described by $\nu$, the solution obtained can be more or less diffused depending on the regularization strength $\gamma$. Moreover, we cannot guarantee the convergence of Algorithm \ref{alg:1} for a fixed weight $\omega$, and to our knowledge, there is no proof for it as of yet. However, if instead of taking fixed values for both $\gamma$ and $\omega$ we consider, at each step $t$, $\gamma(t), \, \omega(t)$ such that $\gamma(t),\,\omega(t)\to 0$ as $t\to+\infty$ and $\sum_t\omega(t)=+\infty$, we can ensure its convergence \citep{IBP}\citep{PeyreGradientFlows}. We have introduced the second condition on the weight $\omega(t)$ to prevent the parameter from vanishing too quickly. Otherwise, in the computation of the WB, we would obtain the final distribution or one close to it, but the subsequent projections introducing the graph constraints could prevent us from reaching such measure, since we may still have no access to those target nodes.
Despite that, in the simulations carried out in Section \ref{sec:4numerical_simulations}, we consider the weight $\omega$ to be both tending to zero (without vanishing too fast) and fixed, since we have observed how, for a constant $\omega<1/2$, the mass reaches the target distribution as well.

\algrenewcommand\algorithmicindent{1.0em}%
\begin{algorithm}[t]
	\linespread{1.1}\selectfont
	\caption{Conceptual procedure of the proposed approach}\label{alg:1}
	\textbf{Input:} Initial and final distributions $\rho_0$ and $\nu$, adjacency matrix $A$, full cost matrix $C^*$, full vector of storage capacities $\rho^*$, regularization parameter $\gamma(t)$ and weight $\omega(t)$ depending on $t$ and such that $\gamma(t),\,\omega(t)\to 0$ as $t\to+\infty$, accuracy parameter $\varepsilon>0$
	\begin{algorithmic}[1]
		\State $t = 0$
		\While{$\frac{1}{2}\norm{\nu-\rho_t}_1 > \varepsilon$}
		\State Find the support of the new measure $\rho_{t+1}$
		\State Define $C$ as the cost matrix $C^*$ but taking only \hspace*{\algorithmicindent}the columns corresponding to the new support
		\State Define $\rho$ as the vector of storage capacity $\rho^*$ but \hspace*{\algorithmicindent}taking only the elements corresponding to the new \hspace*{\algorithmicindent}support
		\State Define the capacity matrix $\Tilde{C}$ as seen in \eqref{eq:cap_mat}
		\State Compute the WB $\rho_{t+1}$ with weights $\omega_1 = \omega(t)$ and \hspace*{\algorithmicindent}$\omega_2 = 1-\omega(t)$ and the additional support and ca-\hspace*{\algorithmicindent}pacity constraints by using Dykstra’s projection \hspace*{\algorithmicindent}algorithm with initial conditions $\pi_1^{(0)} = \pi_2^{(0)} =$ \hspace*{\algorithmicindent}$e^{-\frac{C}{\gamma(t)}}$ and the proximal operators defined on \eqref{eq:ProxC2}, \hspace*{\algorithmicindent}\eqref{eq:ProxC1} and \eqref{eq:proxC4} (with $\rho$) for both transport plans, \hspace*{\algorithmicindent}and \eqref{eq:proxC3} only for transport plan $\pi_1$ to enforce the \hspace*{\algorithmicindent}capacity constraint $\eqref{eq:prob_st_4}$ with capacity matrix $\Tilde{C}$
		\State $t \gets t+1$
		\EndWhile
	\end{algorithmic}
	\textbf{Output:} $\setb{\rho_t}_t$
\end{algorithm}

We state the following lemma regarding the convergence of the computation of each intermediate distribution in the discrete flow.

\begin{lem}\label{lemma:conv_prox}
	For each step $t$, let $\Tilde{C}$ be the capacity matrix defined in \eqref{eq:cap_mat} such that it verifies $\Tilde{C}^\intercal\boldsymbol{1}>\rho_t$, and let $\rho$ be the retention capacity vector in the constraint set $\,\mathcal{C}_2$ such that $\rho_t < \rho$ (both inequalities are considered element-wise). Then,
	the iterative computation of the proximal steps defined in Propositions \ref{prop:prox_C1}, \ref{prop:prox_C2}, \ref{prop:prox_C3} and \ref{prop:prox_C4} converges to the solution of \eqref{eq:prob_st_1}.
\end{lem}
\begin{pf}
	The condition $\Tilde{C}^\intercal\boldsymbol{1}>\rho_t$ ensures that the mass defined by the initial distribution in the $t$-th step, $\rho_t$, can be moved or even kept still in some of the nodes in its support. Similarly, if $\rho$ verifies $\rho_t < \rho$, then the same initial distribution $\rho_t$ is a feasible solution. In particular, we have $\text{ri}(\mathcal{C}_f) \cap \text{ri}(\mathcal{C}_e) \cap \text{ri}(\mathcal{C}_1) \cap \text{ri}(\mathcal{C}_2)  \cap \text{ri}(\mathcal{C}_3) \not=\varnothing$,
	where $\text{ri}(\mathcal{C})$ is the relative interior of the set $\mathcal{C}$. Thus, by Proposition 3.1 in \citep{PeyreGradientFlows}, the iterative computation of proximal steps converges to the desired solution. \qed
\end{pf}
\begin{rem}
	The conditions on $\Tilde{C}$ and $\rho$ are set only to ensure the existence of a feasible solution. Hence, these hypotheses could be exchanged for other expressions as long as they are not so restrictive that a solution cannot satisfy all the constraints. The ones proposed in the statement of Lemma \ref{lemma:conv_prox} are reasonably lax and could be expected in more practical applications. 
\end{rem}

\section{Numerical simulations}\label{sec:4numerical_simulations}
In this section, we show numerical simulations that provide evidence for the effectiveness of the proposed approach and analyze its performance. We further show how to implement it to solve a supply and demand problem related to a drinking water network (DWN).

\subsection{Implementation details}

Before discussing the results obtained in the simulations, some remarks about the implementation of Algorithm \ref{alg:1} are in order. Algorithm \ref{alg:2} shows the detailed steps of Algorithm \ref{alg:1} to solve Problem \ref{prob:prob_statement}. The computations are carried out using logarithms, as some of the values are of the order of $e^{-1/\gamma}$, so when the regularization is really small we might obtain machine precision issues if we did it outside the logarithmic domain. Moreover, to compute $\ln\left((\pi_iq_j)\boldsymbol{1}\right)$ (the product of matrices is element-wise) using $L_{\pi_i}+L_{q_j} = \ln(\pi_i)+\ln(q_j)$, one
\algrenewcommand\algorithmicindent{0.6em}%
\begin{algorithm}[H] 
	\linespread{1.1}\selectfont
	\caption{Detailed implementation of Algorithm \ref{alg:1}}\label{alg:2}
	\textbf{Input:} Initial and final distributions $\rho_0$ and $\nu$, adjacency matrix $A$, full cost matrix $C^*$, full vector of storage capacities $\rho^*$, regularization parameter $\gamma(t)$ and weight $\omega(t)$ depending on $t$ and such that $\gamma(t),\,\omega(t)\to 0$ as $t\to+\infty$, accuracy parameter $\varepsilon$
	\begin{algorithmic}[1]
		\State $t = 0$
		\While{$\frac{1}{2}\norm{\nu-\rho_t}_1 > \varepsilon$}
		\State $supp_{new} = \textsc{supp}(\Bar{A}\rho_t)$ \label{alg:2line_supp_new}
		\State Define $C$ as the cost matrix $C^*$ but taking only \hspace*{\algorithmicindent}the columns corresponding to $supp_{new}$
		\State Define the capacity matrix $\Tilde{C}$ as seen in \eqref{eq:cap_mat}
		\State $L_{\Tilde{C}} = \ln\Tilde{C}$
		\State $L_{\rho} = \ln\rho$
		\State $L_{\pi_1} = L_{\pi_2} = -\frac{1}{\gamma}C$
		\State $L_{q_1} {=} L_{q_2} {=} L_{q_3} {=} L_{q_4} {=} L_{q_5} {=} L_{q_6} {=} L_{q_7} {=} \ln\left(\boldsymbol{1}\boldsymbol{1}^\intercal\right)$
		\State $L_p = \boldsymbol{1}$
		\State $k = 0$
		\While{$\norm{L_p {-} \ln\left(\pi_1^\intercal\boldsymbol{1}\right)}_1 {>} \varepsilon$ \textbf{or} $\abs{\ln(\boldsymbol{1}^\intercal(\pi_1\boldsymbol{1}))}{>} \varepsilon$} \label{alg:2while2}
		\State $L_{\pi_1}' = L_{\pi_1}$
		\State $L_{\pi_2}' = L_{\pi_2}$
		\If{$k \mod 4 = 0$}
		\State $L_{\pi_1} \gets L_{\pi_1}+L_{q_1} + \left(\ln \rho_t - \ln\left((\pi_1q_1)\boldsymbol{1}\right)\right)\boldsymbol{1}^\intercal$
		\State $L_{q_1} \gets L_{q_1} + (L_{\pi_1}'-L_{\pi_1})$
		\State $L_{\pi_2} \gets L_{\pi_2}+L_{q_2} + \left(\ln\nu - \ln\left((\pi_2q_2)\boldsymbol{1}\right)\right)\boldsymbol{1}^\intercal$
		\State $L_{q_2} \gets L_{q_2} + (L_{\pi_2}'-L_{\pi_2})$
		\ElsIf{$k \mod 4 = 1$}
		\State $L_p \gets \omega\ln\left((\pi_1q_3)^\intercal\boldsymbol{1}\right) + (1-\omega)\ln\left((\pi_2q_4)^\intercal\boldsymbol{1}\right)$
		\State $L_{\pi_1} \gets L_{\pi_1} + L_{q_3} + \boldsymbol{1}\left(L_p - \ln\left((\pi_1q_3)^\intercal\boldsymbol{1}\right)\right)^\intercal$
		\State $L_{q_3} \gets L_{q_3} + (L_{\pi_1}'-L_{\pi_1})$
		\State $L_{\pi_2} \gets L_{\pi_2} + L_{q_4} + \boldsymbol{1}\left(L_p - \ln\left((\pi_2q_4)^\intercal\boldsymbol{1}\right)\right)^\intercal$
		\State $L_{q_4} \gets L_{q_4} + (L_{\pi_2}'-L_{\pi_2})$
		\ElsIf{$k \mod 4 = 2$}
		\State $L_{\pi_1} \gets \min(L_{\pi_1}+L_{q_5}, L_{\Tilde{C}})$
		\State $L_{q_5} \gets L_{q_5} + (L_{\pi_1}'-L_{\pi_1})$
		\Else
		\State $L_{\pi_1} \gets L_{\pi_1} + L_{q_6} + \boldsymbol{1}\left(\min\left(L_\rho - \ln\left((\pi_1q_6)^\intercal\boldsymbol{1}\right) \right)\right)^\intercal$
		\State $L_{q_6} \gets L_{q_6} + (L_{\pi_1}'-L_{\pi_1})$
		\State $L_{\pi_2} \gets L_{\pi_2} + L_{q_7} + \boldsymbol{1}\left(\min\left(L_\rho - \ln\left((\pi_2q_7)^\intercal\boldsymbol{1}\right) \right)\right)^\intercal$
		\State $L_{q_7} \gets L_{q_7} + (L_{\pi_2}'-L_{\pi_2})$
		\EndIf
		\State $k \gets k + 1$
		\EndWhile
		\State $\rho_{t+1} = \exp(L_p)$
		\State Rewrite $\rho_{t+1}$ so that it is an $n$-dimensional vector \hspace*{\algorithmicindent}of all zeros except on $supp_{new}$
		\State $t \gets t+1$
		\EndWhile
	\end{algorithmic}
	\textbf{Output:} $\setb{\rho_t}_t$
\end{algorithm}
can take advantage of the identity $\ln\sum_{i=0}^N a_i = \ln a_0 + \ln\left(1 + \sum_{i=1}^N e^{\ln a_i - \ln a_0}\right),$
where $a_0 \geq a_1 \geq \ldots \geq a_N$.

Additionally, for the loop condition at line \ref{alg:2while2} of Algorithm \ref{alg:2}, we have added the second condition $\abs{\ln(\boldsymbol{1}^\intercal(\pi_1\boldsymbol{1}))}> \varepsilon$, to check if a capacity constraint has been enforced on any position on the transport plan $\pi_1$. This is done to avoid numerical issues where, depending on the precision parameter $\varepsilon$, the first \textit{while} condition might not be verified but the solution has not yet converged to $\mathcal{C}_1$ in \eqref{eq:prob_st_4}.

\subsection{Synthetic examples}

To illustrate the steps described in Algorithm \ref{alg:1}, in Figure \ref{fig:Example_star}, we show a simple example, where we start with a Dirac measure at the center of the graph, whose mass has to be distributed among the outermost nodes. Each subsequent plot shows the intermediate measure obtained after one iteration until the final distribution is reached.

In Figure \ref{fig:cost_star}, on the left we plot the total variation distance between the intermediate distribution $\rho_t$ and the target measure $\nu$, for $\omega(t)$ tending to zero at different rates and also fixed at $\omega(t)=0.1$. In any case, we see how we eventually converge to the final distribution. Due to the symmetrical nature of the network and the probability measures, we observe how for $\omega=0.1$, since it gives more weight to minimizing the distance to $\nu$ rather than the previous distribution, the mass advances until it eventually covers the target in a single step. Similarly, for $\omega(t)=1/\ln t$, the weight decreases at a slow rate, and so the mass is transported gradually until $\omega$ is small enough to cover $\nu$ in a single step. For $\omega(t)=1/t$, the decrease rate is faster, but when it finally starts covering $\nu$, it does so fractionally in a couple of steps, since it is still large enough to give some significant weight to the previous distribution. On the right of Figure \ref{fig:cost_star}, we have the cost of transportation (in other words, the Wasserstein distance) of each step, and we observe how the cost adds up to be similar for each case, and we can reach the same conclusions we had with the study of the total variation distance. In particular, we notice how for $\omega(t)=1/\ln t$, the mass does not move until $\omega$ is small enough at the sixth iteration. From there, the transport is similar to what we have for the other cases. As a side note, this is a suitable illustration of the information that the Wasserstein metric can provide with regard to the difference between two measures in the particular domain they are in, which can be lost when using other metrics.

\placefigure[t]{Figure2}

In this and the following examples, we take $\gamma = 10^{-3}$ (except stated otherwise) to have low diffusion. We could take $\gamma(t)\to0$ as $t\to+\infty$, as we have commented earlier. However, the computational speed significantly decreases as the regularization tends to zero, and with this small fixed value, the results obtained have been satisfactory in terms of precision and convergence speed.

\placefigure{Figure3}
Moreover, as we have just commented, apart from taking $\omega(t)\to0$ as $t\to+\infty$, we have also considered a constant weight $\omega = 0.1$, in favor of the final distribution. This constant weight parameter plays an important role in how the mass is transported along with the graph, since being closer to the final distribution rather than the previous measure in terms of the Wasserstein distance does not mean that once the mass has moved from $\rho_0$ to $\nu$, the followed path is the cheapest. To illustrate this case, Figure \ref{fig:Example_physical_capacities} shows an instance where there are two paths to reach the same node from a certain position, and the whole mass must be sent from one place to the other. Taking $\omega = 0.1$, the left plot shows one iteration without adding extra capacity constraints, which results in the mass being transported through the straight line, as one might surmise. The right plot shows the same setting but with a capacity bound of $0.5$ at each link. This value prevents the mass from moving directly to the closest node, and instead, it is forced to be divided and sent through more than one path. With $\omega = 0.1$, the mass tends much more towards the final distribution than the initial one at each step. This amount of $\omega$ forces the mass to be sent through the available paths, which sets the obtained distribution as close as possible to $\nu$ while verifying the constraints of Problem \eqref{eq:prob_st}. However, if $\omega$ is increased, giving preference to the initial measure, what we would observe in the first iteration is a certain amount of mass being sent through the straight path and the rest staying in place in the initial position since it is closer than the secondary path. Thus, we move everything through the straight path, taking more steps to reach the final destination, rather than using all the paths at our disposal to finish in fewer steps, which is cheaper. This scenario highlights the potential use of this weight parameter to model the relationship between distributions at play and determine how we move through the graph.

\placefigure[t]{Figure4}

\subsection{Real case application: drinking water network}
\subsubsection{Small case study}\label{sec:smallDWN}

Now we proceed to study the case of a DWN. Figure \ref{fig:DWN_small} depicts a basic topology of a generic drinking water transport network. The interaction along the most relevant constitutive elements is described by the water supply from the sources towards the network through pumps or valves, depending on the nature of the particular source (either superficial or underground). Therefore, drinking water is moved using manipulated actuators to fill retention tanks and supply water to demand sectors (city neighborhoods). The reader is referred to \citep{Ocampo-Martinez2013} for further details about this system. Here, this case study is used to discuss and analyze how the proposed approach works and how different parameters can be modified, showing the consequent effects over the whole performance of the considered system.

We note that transporting water through a pipe requiring a pump adds a cost of operation to that edge. This added value can be modeled by including the extra expense into the cost matrix.

Figure \ref{fig:Example_small_DWN} shows a simulation on the small network in Figure \ref{fig:DWN_small}, ignoring the pumps (so, no additional costs on the edges). Here, we again take advantage of the parameter $\omega$ to regulate how the water is transported. In particular, in the first step, we use a fairly high weight $\omega = 0.75$ in favor of the initial distribution so that the transportation is done more gradually. In the following steps, as each one is independent of the preceding iteration, the weight is reduced to $\omega=0.1$ so that the demand is covered much faster. Similarly to Figure \ref{fig:cost_star}, Figure \ref{fig:tv_dist_small_DWN} shows the total variation and Wasserstein distance between $\rho_t$ and $\nu$ at each iteration $t$, and we see how we eventually converge to the final distribution with the different weight functions $\omega(t)$ considered. In this case with the chosen network topology and distributions, for $\omega(t)=0.1$ and $\omega(t) = 1/t$, the transportation is identical.

\placefigure[t]{Figure5}
\placefigure[t]{Figure6}

\placefigure[t]{Figure7}

We have seen that with Algorithm \ref{alg:1} presented as it is, we can account for some additional constraints regarding physical limitations, such as capacities on the pipes or additional costs to operate pumps to be able to send mass between certain locations. Constraint \eqref{eq:prob_st_5} in particular has been added with DWN-modelling in mind. According to \citep{Ocampo-Martinez2013}, the nodes that are neither tanks nor sources cannot hold as much water, but they do have a certain retention capacity.

Another issue one can find is having to update specific parameters due to external factors, for instance, the initial or final distributions if, for example, there is a sudden peak in demand, or even the graph topology if a pipe breaks or needs to be cut for maintenance. In the former case, since the scheme is memoryless, the initial and final distributions at any step can be changed, and the algorithm will proceed from there without having to make any modifications to it. For the latter, a change in the topology means that the adjacency and cost matrices are updated, so, as long as these updated values are provided at that step, just as with the change of distributions, the algorithm automatically adapts and proceeds with the computations since the support and capacity constraints are computed at each iteration.

This last case highlights this feature in our approach that we have mentioned several times: each step does not depend on the previous one, which allows the algorithm to adapt to different changes as it advances. If, for example, we wanted to find our sequence of distributions $\setb{\rho_t}_{t\geq0}$ by solving a minimum-cost flow problem, since the flow is computed all at once, each change in the middle of the transportation would mean having to recompute the whole solution (or at least restart taking as the initial measure the distribution obtained at that stage). Simultaneously, with our approach, we only need to update the affected parameters, and the algorithm proceeds from there. Here lies the main difference between our computation of a discrete flow and the continuous flow one would obtain by solving a minimum-cost flow problem.

\subsubsection{Performance assessment with the Barcelona drinking water network}

To show the effectiveness of the proposed approach, a bigger version of a DWN, particularly the one corresponding to Barcelona (Spain) and its metropolitan area, is considered. In this DWN, the water sources are the Ter and Llobregat rivers regulated at their head by some dams with an overall capacity of 600 cubic hectometres. With four drinking water treatment plants, water from rivers and underground sources (wells) is turned into potable water and served to Barcelona and surrounding towns. Those different water sources currently provide a raw flow of around $7$ m$^3$/s. Water flow from each source is limited, implying different water prices depending on water treatments and legal extraction canons.

The Barcelona DWN is structurally organized in two functional layers: an upper layer named \textit{transport network} links the water treatment plants with the reservoirs distributed all over the city, while a lower layer, named \textit{distribution network}, links a reservoir with each consumer sector (water demand). Notice that the upper layer can be managed using control approaches, while the distribution system follows a pre-established behavior given by the water pressure determined. Figure \ref{fig:diagram63tanks} depicts the whole scheme of the transport network. 

Our objective is to implement our algorithm for the management of the upper layer. The setting is analogous to what we have seen in Section \ref{sec:smallDWN} for the small case study: we want to find the (discrete) flow that moves the mass from an initial distribution (water provided by the treatment plants and reservoirs) to a target distribution (expected water in the reservoirs to cover the consumers' water demand) such that it follows the sparsity pattern and constraints induced by the network, and each step is the most cost-efficient (depending on the weight parameter $\omega$). By computing the discrete flow, we can also adapt the solution's next step to changes on the network or the other agents.

To perform the simulations, for the initial distribution $\rho_0$ we have taken the set of source nodes together with close to half of the total amount of tanks (selected at random), assigned them a value following a uniform distribution, and normalized the obtained vector so that $\rho_0\in \text{Prob}(V)$. The final distribution $\nu$ is computed following the same steps with the remaining tanks. For the nodes that are neither tanks nor sources, we have considered that those on the periphery have a retention capacity of $0.05$. For the weight parameter, we have tested it first with a small value $\omega = 0.1$ so that the final distribution is reached in fewer iterations, and then with a larger value $\omega=0.45$, so that the transport is slightly more gradual. Further below we also comment on the convergence when taking $\omega(t)=1/t$ and $\omega(t)=1/\ln t$.

For comparison, the sequence $\setb{\rho_t}_{t\geq0}$ is found by solving Problem \ref{prob:prob_statement}, on one side with Algorithm \ref{alg:1}, using different values of the regularization parameter $\gamma$, and on the other, using the CPLEX solver, which uses the dual simplex algorithm with the default parameters (\texttt{MaxIter} $= 9.2234\times 10^{18}$, \texttt{TolFun} $= 10^{-6}$).

Figure \ref{fig:method_comparison} shows on the top plot the total variation distance between the final distribution $\nu$ and the distribution obtained at every iteration with each method. We notice how with low regularization, the solution obtained is really close (in terms of the total variation distance) to the non-regularized solution obtained with CPLEX, as expected, but even with higher values of the regularization parameter ($\gamma = 1,\,10$), there are no noticeable differences, especially in the case with $\omega = 0.1$. However, with higher values ($\gamma=100$), even though the first iterations are close to the other results, the solution eventually becomes too diffused and is not valid in the setting of DWN. The bottom plot shows the running time of each iteration, i.e., the time elapsed to solve Problem \eqref{eq:prob_st} with the new distribution found in the previous step. As expected, the speed of convergence rapidly decreases as $\gamma \to 0$, which is a known issue with this kind of algorithms \citep{solomon2}. Nonetheless, having seen how with higher regularization, the results obtained are really close even to the CPLEX output, it would be safe to consider a small enough constant $\gamma$ instead of taking $\gamma(t)\to0$ as we do in Algorithm \ref{alg:1}, in exchange of higher performance speed and without losing too much accuracy.
\placefigure{Figure8}

Moreover, we have noticed how by removing the capacity constraint to enforce both \eqref{eq:prob_st_4} and \eqref{eq:prob_st_6}, the algorithm performance vastly improves in terms of convergence speed, which makes sense, considering that it can force sharp changes on the transport plan. In this regard, it would be interesting to find a different approach to improve the computation of the projection with the capacity matrix $\Tilde{C}$ in \eqref{eq:cap_mat}, or directly bypass it by rethinking the constraint in terms of the other parameters and variables at play.

Figure \ref{fig:BCN_steps} shows some selected iterations illustrating how the water is transported towards the target distribution (the bigger the point, the higher the amount of resource is held in that node). Figure \ref{fig:convergence_graphs} shows the total variation distance between $\rho_t$ and $\nu$ at each iteration $t$, taking $\omega(t)=1/t$ (left) and $\omega(t) = 1/\ln t$ (right). As one might expect, since for $\omega(t) = 1/t$ the weight tends to zero at a higher rate, we reach the solution in fewer iterations than taking $\omega(t) = 1/\ln t$. Since the Barcelona DWN is highly connected to cover the whole city and metropolitan area and accounts for any incidents on the network, we have also carried out simulations in different graphs of similar dimensions (around $10^2$ nodes), shown in Figure \ref{fig:convergence_graphs} for comparison. In any case, we observe how the total variation distance eventually converges to zero, taking more steps for the case where the weight decreases slower ($\omega(t) = 1/\ln t$).

\placefigurefull[t]{Figure9}

\placefigurefull{Figure10}

\placefigure[t]{Figure11}

From the point of view related to the management of a DWN, in particular, the considered case of Barcelona, the proposed approach opens new ways of improving existent management criteria in the sense of scalability and modularity of the control approaches 
\citep{Tedesco2018}, apart from adding robustness capabilities to the system. This latter aspect has been previously reported for the particular case given the importance of rejecting the system disturbances and their nature (water costumers demands) 
\citep{Grosso2017}. In any case, a straightforward comparison with existing methods for management and control of DWNs is nowadays not fair since our approach is presented as a proof of concept for the proposed objectives related to the case study, and then some additional design criteria should be considered.

\section{Concluding remarks and future work}\label{sec:5conclusions}

In this paper, we have presented a mathematical formulation to resolve discrete optimal flows over networks based on the computation of constrained Wasserstein Barycenters. Using the entropically regularized approximation of the Wasserstein metric allows us to use Dykstra's projection algorithm, which is easy to implement and is competitive in terms of performance speed since it only requires elementary operations such as matrix and vector products. Moreover, with this methodology, the solution obtained is unique.

We have observed how modifying the capacity matrix to avoid sending mass between non-neighboring nodes forces sharp changes on the transport plan
entries, drastically decreasing the execution speed of the algorithm. Future work should be finding an efficient approach to ensure that this condition is verified. However, this paper focuses on the application of these optimal transport concepts in the context of more \textit{real-life} scenarios and how they can automatically adapt to sudden changes in the topology of the networks or the parameters and distributions.

We have illustrated how the value of the weight $\omega$ alters how the mass is transported from node to node, even mimicking the behavior we can observe if we implement additional physical capacities on the links. It would be interesting to gain more insight into the weight parameter's role in shaping the resulting distribution, not only in our setting but also in the multi-marginal case, with several weights. Moreover, the fact that the methodology proposed can be extended to consider more than two distributions and can adapt to different changes could be used to tackle problems involving decentralized or distributed models, where not all the information is available for every agent.

\begin{dfigurefull}{Figure1}
	\centering
	\subfigure[Initial distribution]{
		\includegraphics[width=0.2\hsize,trim={0cm 0cm 0cm 0.75cm},clip]{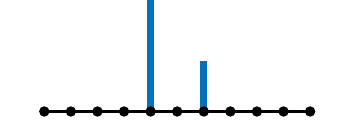}
	}
	\subfigure[Iteration 1]{
		\includegraphics[width=0.2\hsize,trim={0cm 0cm 0cm 0.75cm},clip]{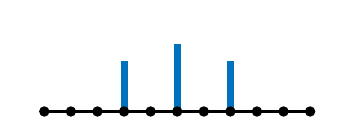}
	}
	\subfigure[Iteration 2]{
		\includegraphics[width=0.2\hsize,trim={0cm 0cm 0cm 0.75cm},clip]{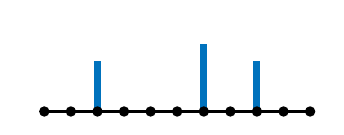}
	}
	\subfigure[Iteration 3]{
		\includegraphics[width=0.2\hsize,trim={0cm 0cm 0cm 0.75cm},clip]{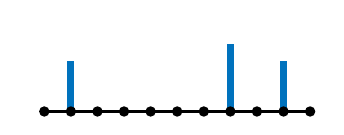}
	}
	\subfigure[Iteration 4]{
		\includegraphics[width=0.2\hsize,trim={0cm 0cm 0cm 0.5cm},clip]{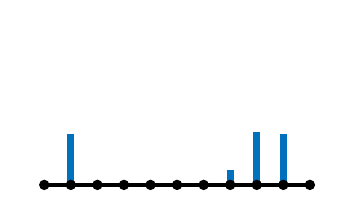}
	}
	\subfigure[Iteration 5]{
		\includegraphics[width=0.2\hsize,trim={0cm 0cm 0cm 0.5cm},clip]{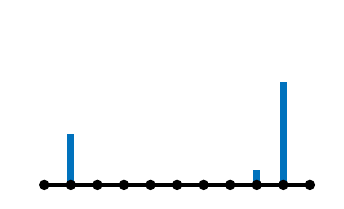}
	}
	\subfigure[Iteration 6]{
		\includegraphics[width=0.2\hsize,trim={0cm 0cm 0cm 0.5cm},clip]{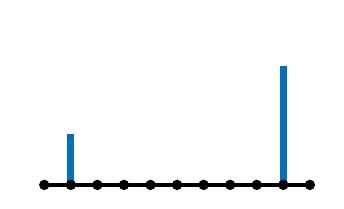}
	}
	\subfigure[Final distribution]{
		\includegraphics[width=0.2\hsize,trim={0cm 0cm 0cm 0.5cm},clip]{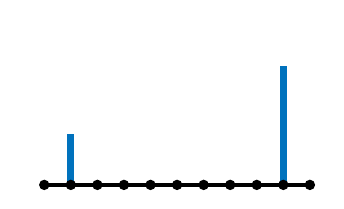}
	}
	\caption{Illustrative example showing the steps obtained by solving problem \ref{prob:prob_statement}. The top left plot shows the initial distribution, and on the bottom right, we have the target distribution, which is reached at iteration 6. Note that there is a capacity constraint on the third-to-last node.}
	\label{fig:example_prob_st}
\end{dfigurefull}


\begin{dfigure}{Figure2}
	\centering
	\subfigure[Initial distribution]{
		\includegraphics[width=0.45\hsize,trim={1.5cm 0.6cm 1.2cm 1.1cm},clip]{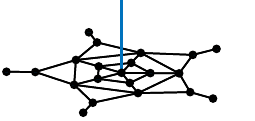}
	}
	\subfigure[Iteration 1]{
		\includegraphics[width=0.45\hsize,trim={1.5cm 0.6cm 1.2cm 1cm},clip]{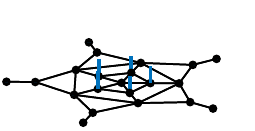}
	}
	\subfigure[Iteration 2]{
		\includegraphics[width=0.45\hsize,trim={1.5cm 0.6cm 1.2cm 1cm},clip]{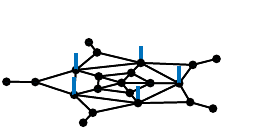}
	}
	\subfigure[Iteration 3]{
		\includegraphics[width=0.45\hsize,trim={1.5cm 0.6cm 1.2cm 1cm},clip]{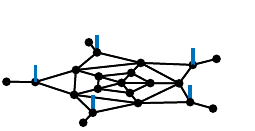}
	}
	\subfigure[Iteration 4]{
		\includegraphics[width=0.45\hsize,trim={1.5cm 0.6cm 1.2cm 1cm},clip]{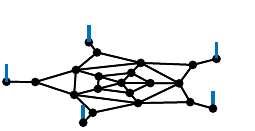}
	}
	\subfigure[Final distribution]{
		\includegraphics[width=0.45\hsize,trim={1.5cm 0.6cm 1.2cm 1cm},clip]{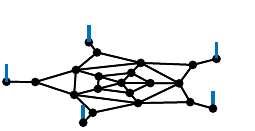}
	}
	\caption{Simple example showing the steps obtained with the algorithm to go from the (top left) initial distribution to the (bottom right) final one.}
	\label{fig:Example_star}
\end{dfigure}


\begin{dfigure}{Figure3}
	\centering
	\begin{tikzpicture}
		\tikzstyle{every node}=[font=\small]
		\begin{axis}[
			width=0.52\hsize,height=3.5cm,scale=1,
			xlabel={Iteration $t$},
			ylabel={$\frac{1}{2}\norm{\nu-\rho_t}_1$},
			ticklabel style = {font=\tiny},
			ymin=0 , ymax=1,
			ytick={0,0.5,1},
			ylabel near ticks,
			xlabel near ticks,
			xmin = 1, xmax=10, 
			every axis plot/.append style={line width=1.3pt}],
			legend pos=south west;
			\addplot  [Black]         	table [x index=0,y index=1]{star_TV.dat};
			\addplot  [NavyBlue]         	table [x index=0,y index=2]{star_TV.dat};
			\addplot  [Orange]       	table [x index=0,y index=3]{star_TV.dat};
		\end{axis}
		\tikzstyle{every node}=[font=\small]
		\begin{axis}[
			legend style={at={(-0.2,-0.5)},anchor=north,nodes={scale=0.65, transform shape}},
			legend columns=-1,
			legend cell align={left},
			width=0.52\hsize,height=3.5cm,scale=1,
			xshift=+4cm,
			xlabel={Iteration $t$},
			ylabel={$\mathcal{W}_\gamma(\rho_{t-1},\rho_t)$},
			ymin=0 , ymax=20,
			ytick={0,10,20},
			xlabel near ticks,
			ylabel near ticks,
			ticklabel style = {font=\tiny},
			xmin = 1, xmax=10, 
			every axis plot/.append style={line width=1.3pt}],
			\addplot  [Black]         	table [x index=0,y index=1]{star_cost.dat};
			\addlegendentry{$\omega(t)=0.1\quad$}
			\addplot  [NavyBlue]         	table [x index=0,y index=2]{star_cost.dat};
			\addlegendentry{$\omega(t)=1/t\quad$}
			\addplot  [Orange]       	table [x index=0,y index=3]{star_cost.dat};
			\addlegendentry{$\omega(t)=1/\ln t$}
		\end{axis}
	\end{tikzpicture}
	\caption{(Left) Total variation distance between the distribution obtained at iteration $t$ ($\rho_t$) and the final distribution ($\nu$), and (right) cost of transportation for each iteration, for the example depicted in Figure \ref{fig:Example_star} and taking different weight functions $\omega(t)$.}
	\label{fig:cost_star}
\end{dfigure}


\begin{dfigure}{Figure4}
	\centering
	\subfigure[(a)]{
		\includegraphics[width=0.42\hsize,trim={0.5cm 0.8cm 1cm 0.35cm},clip]{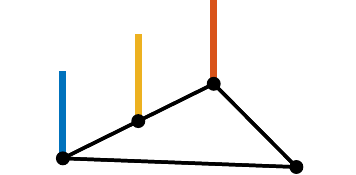}
	}
	\subfigure[(b)]{
		\includegraphics[width=0.42\hsize,trim={0.5cm 0.8cm 1cm 0.35cm},clip]{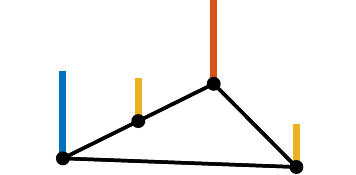}
	}
	\caption{(a) First step computed without adding physical capacities to the links, and (b) with an additional capacity of $0.5$ at each link. We have the initial, final and intermediate distributions in blue, orange, and yellow, respectively.}
	\label{fig:Example_physical_capacities}
\end{dfigure}

\begin{dfigure}{Figure5}
	\centering
	\includegraphics[scale=0.5]{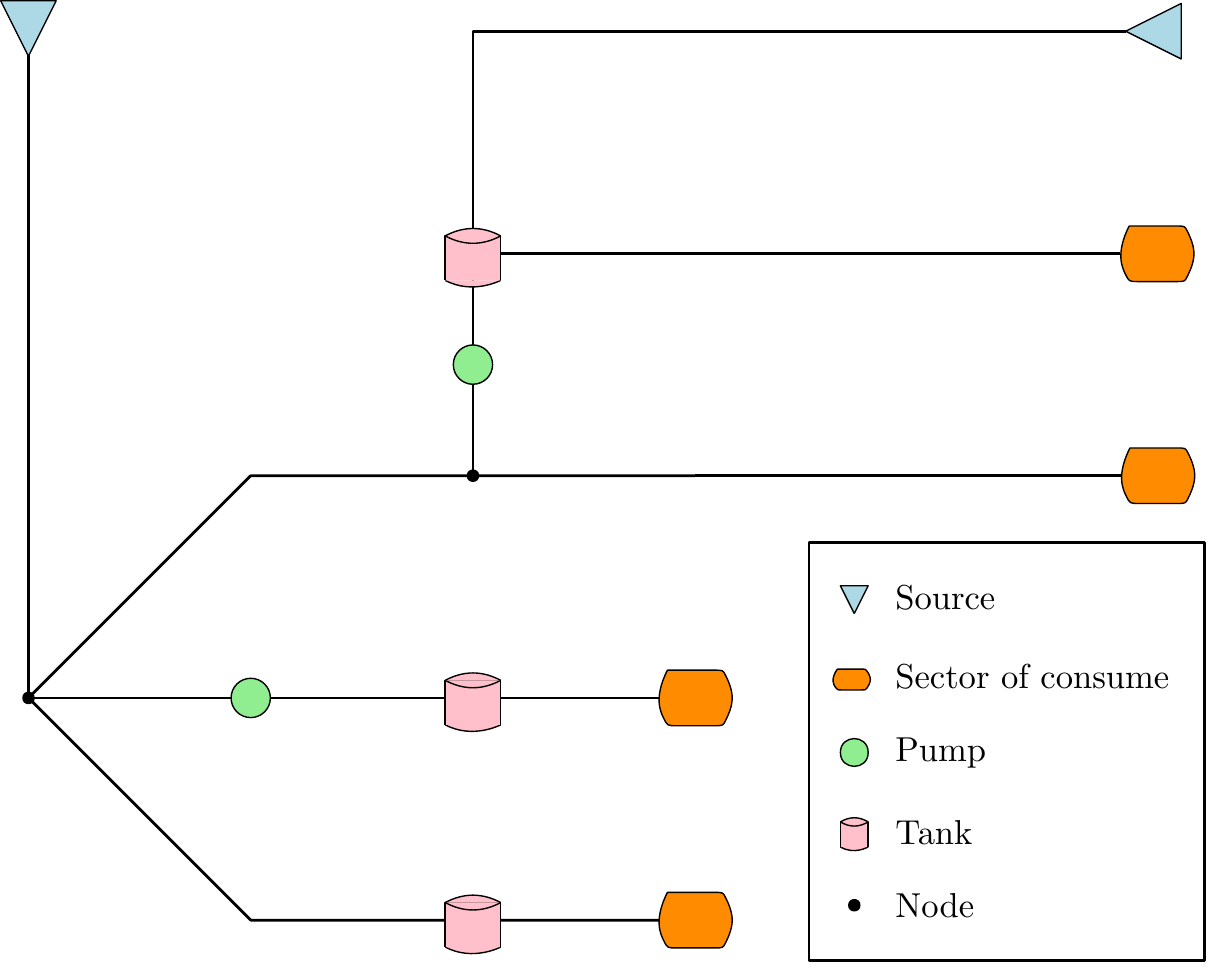}
	\caption{Topology of the small DWN case study.}
	\label{fig:DWN_small}
\end{dfigure}

\begin{dfigure}{Figure6}
	\centering
	\subfigure[Initial distribution]{
		\includegraphics[width=0.45\hsize,trim={0.8cm 0.5cm 0.8cm 1cm},clip]{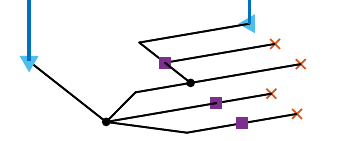}
	}
	\subfigure[Iteration 1, $\omega =0.75$ ]{
		\includegraphics[width=0.45\hsize,trim={0.8cm 0.5cm 0.8cm 1cm},clip]{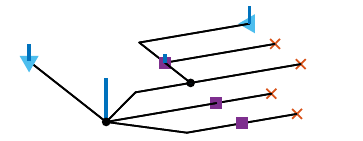}
	}
	\subfigure[Iteration 2, $\omega =0.1$ ]{
		\includegraphics[width=0.45\hsize,trim={0.8cm 0.5cm 0.8cm 1cm},clip]{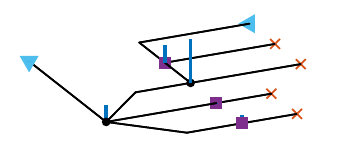}
	}
	\subfigure[Iteration 3, $\omega =0.1$ ]{
		\includegraphics[width=0.45\hsize,trim={0.8cm 0.5cm 0.8cm 1cm},clip]{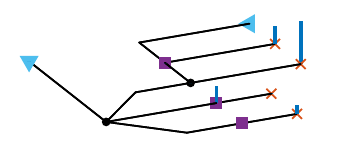}
	}
	\subfigure[Iteration 4, $\omega =0.1$ ]{
		\includegraphics[width=0.45\hsize,trim={0.8cm 0.5cm 0.8cm 1cm},clip]{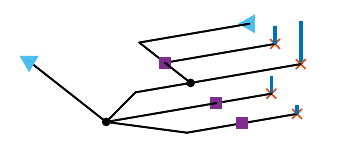}
	}
	\subfigure[Final distribution]{
		\includegraphics[width=0.45\hsize,trim={0.8cm 0.5cm 0.8cm 1cm},clip]{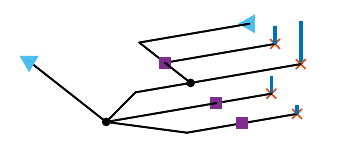}
	}
	\caption{Steps obtained for the small DWN case study.}
	\label{fig:Example_small_DWN}
\end{dfigure}


\begin{dfigure}{Figure7}
	\centering
	\begin{tikzpicture}
		\tikzstyle{every node}=[font=\small]
		\begin{axis}[
			width=0.52\hsize,height=3.4cm,scale=1,
			xlabel={Iteration $t$},
			ylabel={$\frac{1}{2}\norm{\nu-\rho_t}_1$},
			ticklabel style = {font=\tiny},
			ymin=0 , ymax=1,
			ytick={0,0.5,1},
			ylabel near ticks,
			xlabel near ticks,
			xmin=1, xmax=7,
			every axis plot/.append style={line width=1.3pt}],
			legend pos=south west;
			\addplot  [Black]         	table [x index=0,y index=1]{small_DWN_TV.dat};
			\addplot  [dashed,NavyBlue]         	table [x index=0,y index=2]{small_DWN_TV.dat};
			\addplot  [Orange]       	table [x index=0,y index=3]{small_DWN_TV.dat};
		\end{axis}
		\tikzstyle{every node}=[font=\small]
		\begin{axis}[
			legend style={at={(-0.2,-0.5)},anchor=north,nodes={scale=0.65, transform shape}},
			legend columns=-1,
			legend cell align={left},
			width=0.52\hsize,height=3.4cm,scale=1,
			xshift=+4cm,
			xlabel={Iteration $t$},
			ylabel={$\mathcal{W}_\gamma(\rho_{t-1},\rho_t)$},
			ymin=0 , ymax=10,
			ytick={0,5,10},
			xlabel near ticks,
			ylabel near ticks,
			ticklabel style = {font=\tiny},
			xmin=1, xmax=7,
			every axis plot/.append style={line width=1.3pt}],
			\addplot  [Black]         	table [x index=0,y index=1]{small_DWN_cost.dat};
			\addlegendentry{$\omega(t)=0.1\quad$}
			\addplot  [dashed,NavyBlue]         	table [x index=0,y index=2]{small_DWN_cost.dat};
			\addlegendentry{$\omega(t)=1/t\quad$}
			\addplot  [Orange]       	table [x index=0,y index=3]{small_DWN_cost.dat};
			\addlegendentry{$\omega(t)=1/\ln t$}
		\end{axis}
	\end{tikzpicture}
	\caption{(Left) Total variation distance between the distribution obtained at iteration $t$ ($\rho_t$) and the final distribution ($\nu$), and (right) cost of transportation for each iteration, for the example depicted in Figure \ref{fig:Example_small_DWN} (taking different weight functions).}
	\label{fig:tv_dist_small_DWN}
\end{dfigure}


\begin{dfigure}{Figure8}
	\centering
	\begin{tikzpicture}
		\tikzstyle{every node}=[font=\small]
		\begin{axis}[
			width=0.52\hsize,height=3.5cm,scale=1,	
			title = {$\omega=0.1$},
			ylabel={$\frac{1}{2}\norm{\nu-\rho_t}_1$},
			ticklabel style = {font=\tiny},
			ymin=0 , ymax=1,
			ytick={0,0.5,1},
			ylabel near ticks,
			xlabel near ticks,
			xmin = 1, xmax=11, 
			every axis plot/.append style={line width=1.3pt}],
			legend pos=south west;
			\addplot  [NavyBlue]        table [x index=0,y index=1]{comparison_plots_tvd_0p1.dat};
			\addplot  [Orange]         	table [x index=0,y index=2]{comparison_plots_tvd_0p1.dat};
			\addplot  [Yellow]       	table [x index=0,y index=3]{comparison_plots_tvd_0p1.dat};
			\addplot  [Purple]       	table [x index=0,y index=4]{comparison_plots_tvd_0p1.dat};
			\addplot  [LimeGreen]       table [x index=0,y index=5]{comparison_plots_tvd_0p1.dat};
			\addplot  [Cerulean]       	table [x index=0,y index=6]{comparison_plots_tvd_0p1.dat};
		\end{axis}
		\tikzstyle{every node}=[font=\small]
		\begin{axis}[
			yshift=-3.2cm,
			width=0.52\hsize,height=3.5cm,scale=1,
			title = {$\omega=0.1$},
			ylabel={$\mathcal{W}_\gamma(\rho_{t-1},\rho_t)$},
			ticklabel style = {font=\tiny},
			ymin=0 , ymax=300,
			ytick={0,150,300},
			ylabel near ticks,
			xlabel near ticks,
			xmin = 1, xmax=11, 
			every axis plot/.append style={line width=1.3pt}],
			legend pos=south west;
			\addplot  [NavyBlue]        table [x index=0,y index=1]{comparison_plots_cost_0p1.dat};
			\addplot  [Orange]         	table [x index=0,y index=2]{comparison_plots_cost_0p1.dat};
			\addplot  [Yellow]       	table [x index=0,y index=3]{comparison_plots_cost_0p1.dat};
			\addplot  [Purple]       	table [x index=0,y index=4]{comparison_plots_cost_0p1.dat};
			\addplot  [LimeGreen]       table [x index=0,y index=5]{comparison_plots_cost_0p1.dat};
			\addplot  [Cerulean]       	table [x index=0,y index=6]{comparison_plots_cost_0p1.dat};
		\end{axis}
		\begin{axis}[
			legend style={nodes={scale=0.7, transform shape}},
			legend pos=north east,
			legend cell align={left},
			width=0.52\hsize,height=3.5cm,scale=1,
			yshift=-6.4cm,      		
			title = {$\omega=0.1$},
			ymin=-2.7 , ymax=3.72,
			ytick={-2.5,0,2.5},
			ylabel={log(run time (s))},
			xlabel={Iteration $t$},
			xlabel near ticks,
			ylabel near ticks,
			ticklabel style = {font=\tiny},
			xmin = 1, xmax=11,
			every axis plot/.append style={line width=1.3pt}],
			\addplot  [NavyBlue]        table [x index=0,y index=1]{comparison_plots_logtime_0p1.dat};
			\addplot  [Orange]         	table [x index=0,y index=2]{comparison_plots_logtime_0p1.dat};
			\addplot  [Yellow]       	table [x index=0,y index=3]{comparison_plots_logtime_0p1.dat};
			\addplot  [Purple]       	table [x index=0,y index=4]{comparison_plots_logtime_0p1.dat};
			\addplot  [LimeGreen]       table [x index=0,y index=5]{comparison_plots_logtime_0p1.dat};
			\addplot  [Cerulean]       	table [x index=0,y index=6]{comparison_plots_logtime_0p1.dat};
		\end{axis}
		\tikzstyle{every node}=[font=\small]
		\begin{axis}[
			legend style={at={(-0.2,-3.9)},anchor=north,nodes={scale=0.65, transform shape}},
			legend columns=3,
			legend cell align={left},
			xshift=+4cm,
			width=0.52\hsize,height=3.5cm,scale=1,  		
			title = {$\omega=0.45$},
			ymin=0 , ymax=1,
			ytick={0,0.5,1},
			xtick={3,6,9,12,15},
			xlabel near ticks,
			ticklabel style = {font=\tiny},
			xmin = 1, xmax=15, 
			every axis plot/.append style={line width=1.3pt}],
			\addplot  [NavyBlue]        table [x index=0,y index=1]{comparison_plots_tvd_0p45.dat};
			\addlegendentry{CPLEX}
			\addplot  [Orange]         	table [x index=0,y index=2]{comparison_plots_tvd_0p45.dat};
			\addlegendentry{$\gamma = 10^{-2}$}
			\addplot  [Yellow]       	table [x index=0,y index=3]{comparison_plots_tvd_0p45.dat};
			\addlegendentry{$\gamma = 10^{-1}$}
			\addplot  [Purple]       	table [x index=0,y index=4]{comparison_plots_tvd_0p45.dat};
			\addlegendentry{$\gamma = 10^{0}\quad$}
			\addplot  [LimeGreen]         	table [x index=0,y index=5]{comparison_plots_tvd_0p45.dat};
			\addlegendentry{$\gamma = 10^{1}\quad$}
			\addplot  [Cerulean]       	table [x index=0,y index=6]{comparison_plots_tvd_0p45.dat};
			\addlegendentry{$\gamma = 10^{2}$}
		\end{axis}
		\tikzstyle{every node}=[font=\small]
		\begin{axis}[
			legend style={nodes={scale=0.7, transform shape}},
			legend pos=north east,
			legend cell align={left},
			width=0.52\hsize,height=3.5cm,scale=1,
			yshift=-3.2cm,
			xshift=+4cm,   		
			title = {$\omega=0.45$},
			ymin=0 , ymax=300,
			ytick={0,150,300},
			xtick={3,6,9,12,15},
			xlabel near ticks,
			ticklabel style = {font=\tiny},
			xmin = 1, xmax=15, 
			every axis plot/.append style={line width=1.3pt}],
			\addplot  [NavyBlue]        table [x index=0,y index=1]{comparison_plots_cost_0p45.dat};
			\addplot  [Orange]         	table [x index=0,y index=2]{comparison_plots_cost_0p45.dat};
			\addplot  [Yellow]       	table [x index=0,y index=3]{comparison_plots_cost_0p45.dat};
			\addplot  [Purple]       	table [x index=0,y index=4]{comparison_plots_cost_0p45.dat};
			\addplot  [LimeGreen]       table [x index=0,y index=5]{comparison_plots_cost_0p45.dat};
			\addplot  [Cerulean]       	table [x index=0,y index=6]{comparison_plots_cost_0p45.dat};
		\end{axis}
		\begin{axis}[
			legend style={nodes={scale=0.7, transform shape}},
			legend pos=north east,
			legend cell align={left},
			width=0.52\hsize,height=3.5cm,scale=1,
			yshift=-6.4cm,
			xshift=+4cm,     		
			title = {$\omega=0.45$},
			ymin=-2.7 , ymax=3.72,
			ytick={-2.5,0,2.5},
			xtick={3,6,9,12,15},
			xlabel={Iteration $t$},
			xlabel near ticks,
			ticklabel style = {font=\tiny},
			xmin = 1, xmax=15, 
			every axis plot/.append style={line width=1.3pt}],
			\addplot  [NavyBlue]        table [x index=0,y index=1]{comparison_plots_logtime_0p45.dat};
			\addplot  [Orange]         	table [x index=0,y index=2]{comparison_plots_logtime_0p45.dat};
			\addplot  [Yellow]       	table [x index=0,y index=3]{comparison_plots_logtime_0p45.dat};
			\addplot  [Purple]       	table [x index=0,y index=4]{comparison_plots_logtime_0p45.dat};
			\addplot  [LimeGreen]       table [x index=0,y index=5]{comparison_plots_logtime_0p45.dat};
			\addplot  [Cerulean]       	table [x index=0,y index=6]{comparison_plots_logtime_0p45.dat};
		\end{axis}
	\end{tikzpicture}
	\caption{Performance comparison between Algorithm \ref{alg:1} (using increasing values of the regularization parameter) and CPLEX, using $\omega = 0.1$ (left column) and $\omega = 0.45$ (right column). (Top) Total variation distance between the final distribution $\nu$ and the distribution obtained at iteration $t$ ($\rho_t$), (middle) cost of transportation for each iteration, and (bottom) time elapsed (in seconds) for each iteration (the plot is in logarithmic scale for visualization purposes).}
	\label{fig:method_comparison}
\end{dfigure}


\begin{dfigurefull}{Figure9}
	\centering
	\includegraphics[width=\hsize]{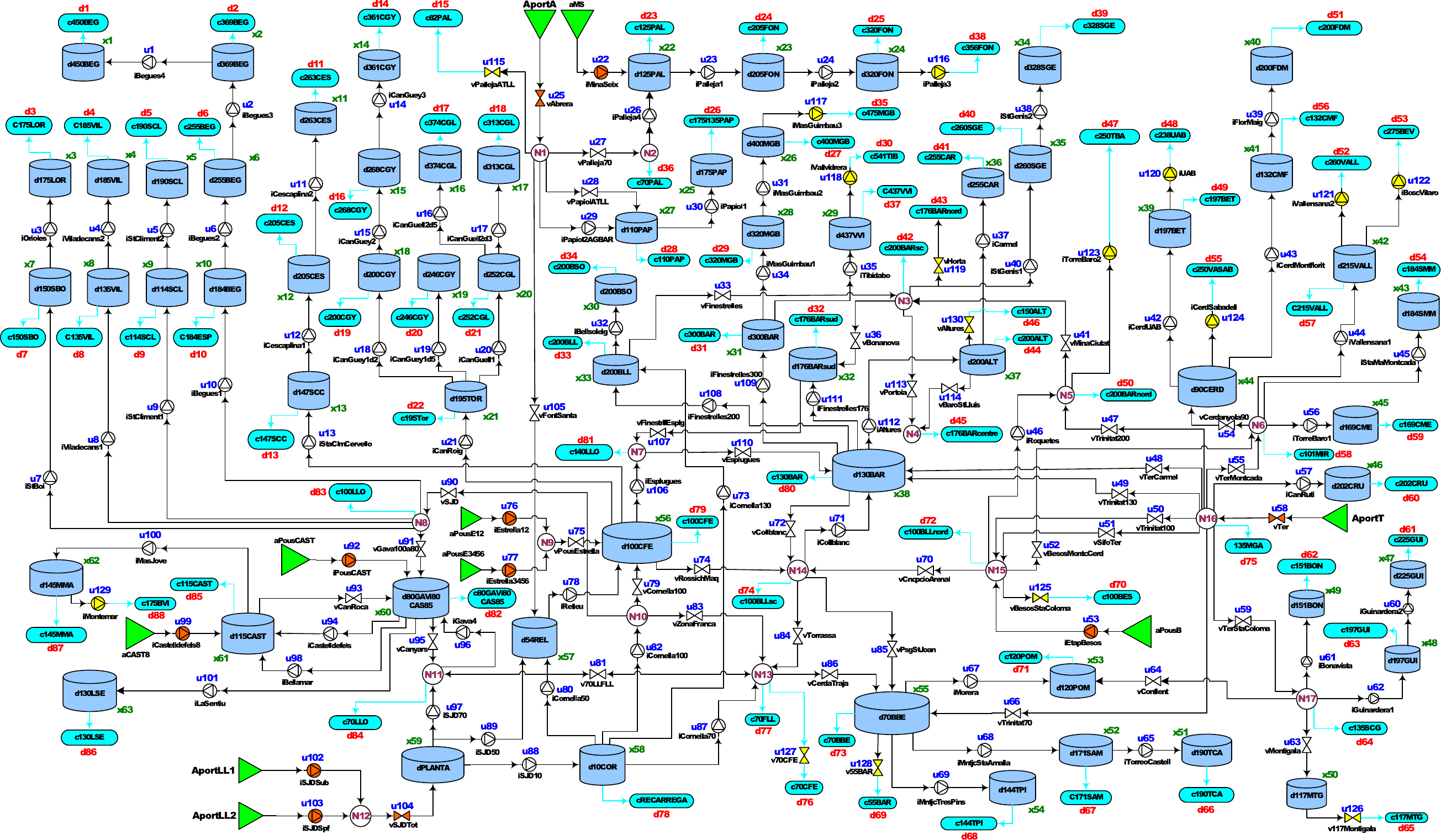}
	\caption{Scheme of the Barcelona drinking water transport network. Taken from \citep{Ocampo-Martinez2013}.}
	\label{fig:diagram63tanks}
\end{dfigurefull}


\begin{dfigurefull}{Figure10}
	\centering
	\subfigure[Initial distribution]{
		\includegraphics[width=0.25\hsize,trim={0cm 0.2cm 0cm 0cm},clip]{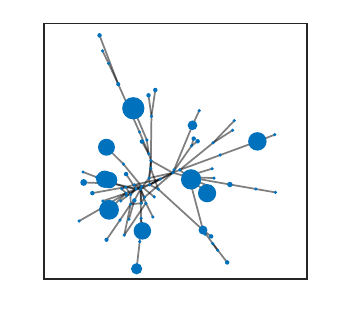}
	}
	\subfigure[Iteration 1]{
		\includegraphics[width=0.25\hsize,trim={0cm 0.2cm 0cm 0cm},clip]{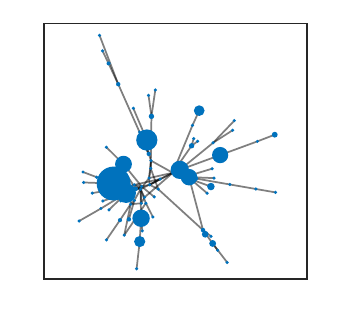}
	}
	\subfigure[Iteration 3]{
		\includegraphics[width=0.25\hsize,trim={0cm 0.2cm 0cm 0cm},clip]{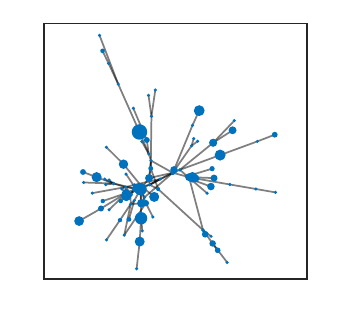}
	}
	\subfigure[Iteration 5]{
		\includegraphics[width=0.25\hsize,trim={0cm 0.2cm 0cm 0cm},clip]{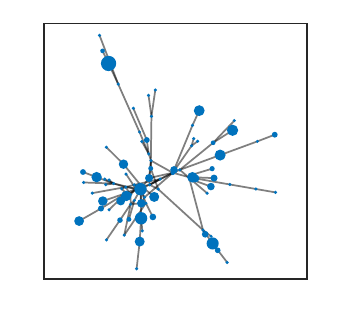}
	}
	\subfigure[Iteration 7]{
		\includegraphics[width=0.25\hsize,trim={0cm 0.2cm 0cm 0cm},clip]{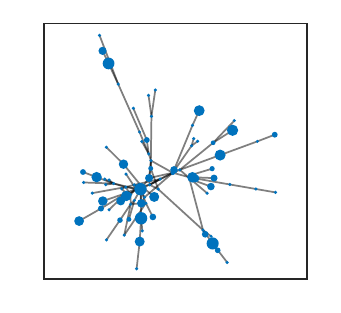}
	}
	\subfigure[Final distribution]{
		\includegraphics[width=0.25\hsize,trim={0cm 0.2cm 0cm 0cm},clip]{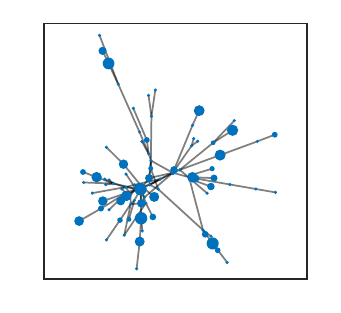}
	}
	\caption{Selected steps showing the transportation of mass through the Barcelona drinking water network from the (top left) initial distribution to the (bottom right) final distribution. Parameters used: $\omega=0.1$ and $\gamma=10^{-1}$.}
	\label{fig:BCN_steps}
\end{dfigurefull}


\begin{dfigure}{Figure11}
	\centering
	\begin{tikzpicture}
		\tikzstyle{every node}=[font=\small]
		\begin{axis}[
			width=0.52\hsize,height=3.5cm,scale=1,	
			title = {$\omega=1/t$},
			xlabel={Iteration $t$},
			ylabel={$\frac{1}{2}\norm{\nu-\rho_t}_1$},
			ticklabel style = {font=\tiny},
			ymin=0 , ymax=1,
			ytick={0,0.5,1},
			ylabel near ticks,
			xlabel near ticks,
			xmin = 1, xmax=20, 
			every axis plot/.append style={line width=1.3pt}],
			legend pos=south west;
			\addplot  [NavyBlue]         	table [x index=0,y index=1]{convergence_plot1.dat};
			\addplot  [Orange]         	table [x index=0,y index=2]{convergence_plot1.dat};
			\addplot  [LimeGreen]       	table [x index=0,y index=3]{convergence_plot1.dat};
			\addplot  [Purple]       	table [x index=0,y index=4]{convergence_plot1.dat};
		\end{axis}
		\tikzstyle{every node}=[font=\small]
		\begin{axis}[
			legend style={at={(-0.2,-0.5)},anchor=north,nodes={scale=0.65, transform shape}},
			legend columns=-1,
			legend cell align={left},
			width=0.52\hsize,height=3.5cm,scale=1,
			xshift=+4cm,		
			title = {$\omega=1/\ln t$},
			xlabel={Iteration $t$},
			ymin=0 , ymax=1,
			ytick={0,0.5,1},
			xlabel near ticks,
			ticklabel style = {font=\tiny},
			xmin = 1, xmax=90, 
			every axis plot/.append style={line width=1.3pt}],
			\addplot  [NavyBlue]         	table [x index=0,y index=1]{convergence_plot2.dat};
			\addlegendentry{BCN DWN$\quad$}
			\addplot  [Orange]         	table [x index=0,y index=2]{convergence_plot2.dat};
			\addlegendentry{Star$\quad$}
			\addplot  [LimeGreen]       	table [x index=0,y index=3]{convergence_plot2.dat};
			\addlegendentry{Grid$\quad$}
			\addplot  [Purple]       	table [x index=0,y index=4]{convergence_plot2.dat};
			\addlegendentry{Cycle}
		\end{axis}
	\end{tikzpicture}
	\caption{Total variation distance between the final distribution ($\nu$) and the one computed at iteration t ($\rho_t$), with $\omega=1/t$ (left) and $\omega=1/\ln t$ (right), computed on different graphs, all of size $10^2$. $\gamma=10^{-1}$ for every case.}
	\label{fig:convergence_graphs}
\end{dfigure}

 \begin{ack}    
This work was partially funded by ARPA-H Strategic Initiative Seed Fund \#916012, Sustainable Futures Fund \#919027, and the Spanish project PID2020-115905RB-C21 (L-BEST) funded by MCIN/ AEI /10.13039/501100011033.
 \end{ack}

\bibliographystyle{agsm}        
\bibliography{autosam}           

@Inbook{user,
    author = {Ambrosio, L. and Gigli, N.},
    title = {A User's Guide to Optimal Transport},
    bookTitle = {Modelling and Optimisation of Flows on Networks: Cetraro, Italy 2009, Editors: Benedetto Piccoli, Michel Rascle},
    year = {2013},
    publisher = {Springer Berlin Heidelberg},
    address = {Berlin, Heidelberg},
    pages = {1--155},
    doi = {10.1007/978-3-642-32160-3_1}
}

@book{gradientFlows,
  title={Gradient Flows: In Metric Spaces and in the Space of Probability Measures},
  author={Ambrosio, L. and Gigli, N. and Savare, G.},
  isbn={9783764387228},
  lccn={2008921489},
  series={Lectures in Mathematics. ETH Z{\"u}rich},
  year={2008},
  publisher={Birkh{\"a}user Basel}
}

@article{GAN,
  title={{W}asserstein {GAN}},
  author={Arjovsky, M. and Chintala, S. and Bottou, L.},
  journal={arXiv},
  year={2017},
  volume={1701.07875}
}

@article{bassetti2,
author = {Bassetti, F. and Gualandi, S. and Veneroni, M.},
title = {On the Computation of Kantorovich--Wasserstein Distances Between Two-Dimensional Histograms by Uncapacitated Minimum Cost Flows},
journal = {SIAM Journal on Optimization},
volume = {30},
number = {3},
pages = {2441-2469},
year = {2020},
doi = {10.1137/19M1261195},
eprint = { 
        https://doi.org/10.1137/19M1261195
    
}

}

@article{brenier,
author = {Brenier, Y.},
title = {Polar factorization and monotone rearrangement of vector-valued functions},
journal = {Communications on Pure and Applied Mathematics},
volume = {44},
number = {4},
pages = {375-417},
doi = {10.1002/cpa.3160440402},
eprint = {https://onlinelibrary.wiley.com/doi/pdf/10.1002/cpa.3160440402},
year = {1991}
}

@inproceedings{cuturiSinkhorn,
  title={Sinkhorn distances: Lightspeed computation of optimal transport},
  author={Cuturi, M.},
  booktitle={Advances in neural information processing systems},
  pages={2292--2300},
  year={2013}
}

@article{cuturiFastWB,
  title={Fast computation of {W}asserstein barycenters},
  author={Cuturi, M. and Doucet, A.},
  year={2014},
  publisher={Journal of Machine Learning Research}
}

@article{solomon2,
  title={Quadratically regularized optimal transport on graphs},
  author={Essid, M. and Solomon, J.},
  journal={SIAM Journal on Scientific Computing},
  volume={40},
  number={4},
  pages={A1961--A1986},
  year={2018},
  publisher={SIAM}
}

@article{kant3,
author = {Kantorovitch, L. V.},
title = {On the Translocation of Masses},
journal = {Management Science},
volume = {5},
number = {1},
pages = {1-4},
year = {1958},
doi = {10.1287/mnsc.5.1.1},
eprint = { 
        https://doi.org/10.1287/mnsc.5.1.1
    
}
}

@article{kovacs,
author = {P. Kovács},
title = {Minimum-cost flow algorithms: an experimental evaluation},
journal = {Optimization Methods and Software},
volume = {30},
number = {1},
pages = {94-127},
year  = {2015},
publisher = {Taylor & Francis},
doi = {10.1080/10556788.2014.895828},
eprint = { 
        https://doi.org/10.1080/10556788.2014.895828
    
}

}

@book{monge,
  title={M{\'e}moire sur la th{\'e}orie des d{\'e}blais et des remblais},
  author={Monge, G.},
  year={1781},
  publisher={De l'Imprimerie Royale}
}

@article{Ocampo-Martinez2013,
	Author = {Ocampo-Martinez, C. and Puig, V. and Cembrano, G. and Quevedo, J.},
	Journal = {IEEE Control Systems Magazine},
	Number = {1},
	Pages = {15-41},
	Title = {Application of Predictive Control Strategies to the Management of Complex Networks in the Urban Water Cycle},
	Volume = {33},
	Year = {2013},
	doi = {10.1109/MCS.2012.2225919}
	}

@article{compOT,
year = {2019},
volume = {11},
journal = {Foundations and Trends in Machine Learning},
title = {Computational Optimal Transport: With Applications to Data Science},
doi = {10.1561/2200000073},
issn = {1935-8237},
number = {5-6},
pages = {355-607},
author = {G. Peyré and M. Cuturi}
}

@article{rubner,
  title={The earth mover's distance as a metric for image retrieval},
  author={Rubner, Y. and Tomasi, C. and Guibas, L. J.},
  journal={International journal of computer vision},
  volume={40},
  number={2},
  pages={99--121},
  year={2000},
  publisher={Springer}
}

@article{solomonConvolutional,
author = {Solomon, J. and de Goes, F. and Peyr\'{e}, G. and Cuturi, M. and Butscher, A. and Nguyen, A. and Du, T. and Guibas, L. J.},
title = {Convolutional {W}asserstein Distances: Efficient Optimal Transportation on Geometric Domains},
year = {2015},
issue_date = {July 2015},
publisher = {Association for Computing Machinery},
address = {New York, NY, USA},
volume = {34},
number = {4},
issn = {0730-0301},
doi = {10.1145/2766963},
journal = {ACM Trans. Graph.},
month = jul,
articleno = {Article 66},
numpages = {11}
}

@book{villani,
  title={Optimal Transport: Old and New},
  author={Villani, C.},
  isbn={9783540710509},
  lccn={2008932183},
  series={Grundlehren der mathematischen Wissenschaften},
  year={2008},
  publisher={Springer Berlin Heidelberg}
}

@article{IBP,
author = {Benamou, J.-D. and Carlier, G. and Cuturi, M. and Nenna, L. and Peyré, G.},
title = {Iterative {B}regman Projections for Regularized Transportation Problems},
journal = {SIAM Journal on Scientific Computing},
volume = {37},
number = {2},
pages = {A1111-A1138},
year = {2015},
doi = {10.1137/141000439},
eprint = { 
        https://doi.org/10.1137/141000439
    
}

}

@article{Dykstra,
 ISSN = {01621459},
 author = {R. L. Dykstra},
 journal = {Journal of the American Statistical Association},
 number = {384},
 pages = {837--842},
 publisher = {[American Statistical Association, Taylor & Francis, Ltd.]},
 title = {An Algorithm for Restricted Least Squares Regression},
 volume = {78},
 year = {1983}
}

@article{Chow2017ADS,
  title={A discrete {S}chrodinger equation via optimal transport on graphs},
  author={S.-N. Chow and W. Li and H.-M. Zhou},
  journal={arXiv: Dynamical Systems},
  year={2017}
}

@article{Richemond2017OnWR,
  title={On {W}asserstein Reinforcement Learning and the {F}okker-{P}lanck equation},
  author={P. H. Richemond and B. Maginnis},
  journal={ArXiv},
  year={2017},
  volume={abs/1712.07185}
}

@article{Mielke2013GeodesicCO,
  title={Geodesic convexity of the relative entropy in reversible {M}arkov chains},
  author={A. Mielke},
  journal={Calculus of Variations and Partial Differential Equations},
  year={2013},
  volume={48},
  pages={1-31}
}

@article{Cuturi2016ASD,
  title={A Smoothed Dual Approach for Variational {W}asserstein Problems},
  author={M. Cuturi and G. Peyr{\'e}},
  journal={SIAM J. Imaging Sci.},
  year={2016},
  volume={9},
  pages={320-343}
}

@article{PeyreGradientFlows,
author = {Peyré, G.},
title = {Entropic Approximation of {W}asserstein Gradient Flows},
journal = {SIAM Journal on Imaging Sciences},
volume = {8},
number = {4},
pages = {2323-2351},
year = {2015},
doi = {10.1137/15M1010087},
eprint = { 
        https://doi.org/10.1137/15M1010087
    
}

}

@Article{Tedesco2018,
	author = {F. Tedesco and C. Ocampo-Martinez and A. Cassavola and V. Puig},
	title = {Centralised and Distributed Command Governor Approaches for the Operational Control of Drinking Water Networks},
	journal = {IEEE Transactions on Systems, Man \& Cybernetics: Systems},
	year = {2018},
	volume = {48},
	number = {4},
	pages = {586-595},
	doi = {10.1109/TSMC.2016.2612361}
}

@Article{Grosso2017,
	author = {J.M. Grosso and P. Velarde and C. Ocampo-Martinez and J.M. Maestre and V. Puig},
	title = {Stochastic model predictive control approaches applied to drinking water networks},
	journal = {Optimal Control, Applications and Methods},
	year = {2017},
	volume = {38},
	number = {4},
	pages = {541-558},
	doi = {10.1002/oca.2269}
}

@article{erbar,
  title={Computation of optimal transport on discrete metric measure spaces},
  author={M. Erbar and M. Rumpf and B. Schmitzer and S. Simon},
  journal={Numerische Mathematik},
  year={2020},
  volume={144},
  pages={157-200}
}

@article{JKO,
author = {J. Richard and K. David and O. Felix},
title = {The Variational Formulation of the {F}okker-{P}lanck Equation},
journal = {SIAM Journal on Mathematical Analysis},
volume = {29},
number = {1},
pages = {1-17},
year = {1998},
doi = {10.1137/S0036141096303359},
eprint = { 
        https://doi.org/10.1137/S0036141096303359
    
}

}

@article{Vaserstein,
  title={Markov processes over denumerable products of spaces describing large systems of automata},
  author={L. N. Vaserstein},
  journal = {Problems of Information Transmission},
  volume = {5},
  number = {3},
  pages = {47–52},
  year={1969}
}

@article{MCF1967,
author = {K. Morton},
title = {A Primal Method for Minimal Cost Flows with Applications to the Assignment and Transportation Problems},
journal = {Management Science},
volume = {14},
number = {3},
pages = {205-220},
year = {1967},
doi = {10.1287/mnsc.14.3.205},
eprint = {https://doi.org/10.1287/mnsc.14.3.205},
}

@book{MCF1993,
  title={Network Flows: Theory, Algorithms, and Applications},
  author={Ahuja, R.K. and Magnanti, T.L. and Orlin, J.B.},
  isbn={9780136175490},
  lccn={lc92026702},
  year={1993},
  publisher={Prentice Hall}
}

@article{MCF2020,
  title={An Efficient Algorithm for Solving Minimum Cost Flow Problem with Complementarity Slack Conditions},
  author={Y. Hu and X. Zhao and J. Liu and B. Liang and C. Ma},
  journal={Mathematical Problems in Engineering},
  year={2020},
  volume={2020},
  pages={1-5}
}

@article{matrixFact,
  title={Supervised Quantile Normalization for Low-rank Matrix Approximation},
  author={M. Cuturi and O. Teboul and J. Niles-Weed and J.-P. Vert},
  journal={arXiv},
  year={2020},
  volume={2002.03229}
}

@article{FairML,
  title={Fair Regression with Wasserstein Barycenters},
  author={E. Chzhen and C. Denis and M. Hebiri and L. Oneto and M. Pontil},
  journal={arXiv},
  year={2020},
  volume={2006.07286}
}

@article{averagingGas,
  title={Averaging Atmospheric Gas Concentration Data using {W}asserstein Barycenters},
  author={M. Barré and C. Giron and M. Mazzolini and A. d'Aspremont},
  journal={arXiv},
  year={2020},
  volume={2010.02762}
}

@inproceedings{neuronLabel,
  title={Probabilistic Joint Segmentation and Labeling of C. elegans Neurons},
  author={Nejatbakhsh, A. and Varol, E. and Yemini, E. and Hobert, O. and Paninski, L.},
  booktitle={International Conference on Medical Image Computing and Computer-Assisted Intervention},
  pages={130--140},
  year={2020},
  organization={Springer}
}

@article{covid,
  title={Clustering patterns connecting COVID-19 dynamics and Human mobility using optimal transport},
  author={F. Nielsen and G. Marti and S. Ray and S. Pyne},
  journal={arXiv},
  year={2020},
  volume={2007.10677}
}

@inproceedings{decentralizeAndRand,
title = {Decentralize and Randomize: Faster Algorithm for Wasserstein Barycenters},
author = {Dvurechenskii, P. and Dvinskikh, D. and Gasnikov, A. and Uribe, C. and Nedich, A.},
booktitle = {Advances in Neural Information Processing Systems 31},
pages = {10760--10770},
year = {2018},
publisher = {Curran Associates, Inc.}
}

@article{distributedQuantization,
  title={Distributed Optimization with Quantization for Computing {W}asserstein Barycenters},
  author={Krawtschenko, R. and Uribe, C. and Gasnikov, A. and Dvurechensky, P.},
  journal={arXiv},
  year={2020},
  volume={2010.14325}
}

@article{multimarginalOT,
  title={Multimarginal Optimal Transport by Accelerated Gradient Descent},
  author={N. Tupitsa and P. Dvurechensky and A. Gasnikov and C. Uribe},
  journal={arXiv},
  year={2020},
  volume={2004.02294}
}

@Article{Bunne2021,
  author        = {C. Bunne and L. Meng-Papaxanthos and . Krause and M. Cuturi},
  title         = {JKOnet: Proximal Optimal Transport Modeling of Population Dynamics},
  year          = {2021},
  journal       = {arXiv},
  volume        = {2106.06345},
}

@Article{Haasler2021,
  author    = {Isabel Haasler and Yongxin Chen and Johan Karlsson},
  journal   = {{IEEE} Control Systems Letters},
  title     = {Optimal Steering of Ensembles With Origin-Destination Constraints},
  year      = {2021},
  issn      = {2475-1456},
  number    = {3},
  pages     = {881-886},
  volume    = {5},
  doi       = {10.1109/lcsys.2020.3006763}
}

@Article{Guex2019,
  author    = {Guillaume Guex and Ilkka Kivim{\"{a}}ki and Marco Saerens},
  journal   = {Netw. Sci.},
  title     = {Randomized optimal transport on a graph: framework and new distance measures},
  year      = {2019},
  month     = {mar},
  number    = {1},
  pages     = {88--122},
  volume    = {7},
  bibsource = {dblp computer science bibliography, https://dblp.org},
  doi       = {10.1017/nws.2018.29},
  publisher = {Cambridge University Press ({CUP})},
}

@Article{Chen2018,
  author    = {Yongxin {Chen} and Tryphon T. {Georgiou} and Michele {Pavon} and Allen {Tannenbaum}},
  journal   = {{IEEE} Transactions on Automatic Control},
  title     = {Efficient robust routing for single commodity network flows},
  year      = {2018},
  issn      = {0018-9286},
  month     = {jul},
  number    = {7},
  pages     = {2287--2294},
  volume    = {63},
  doi       = {10.1109/tac.2017.2763418},
  language  = {English},
  publisher = {Institute of Electrical and Electronics Engineers ({IEEE})},
  zbl       = {1423.90029},
}

@Book{Schroedinger1931,
  author    = {E. {Schr\"odinger}},
  publisher = {Verlag der Akademie der Wissenschaften in Kommission bei Walter De Gruyter u. Co.},
  title     = {Über die Umkehrung der Naturgesetze},
  year      = {1931},
  address   = {Berlin},
  volume    = {9},
  journal   = {Sitzungsberichte der Preu{\ss}ischen Akademie der Wissenschaften, Physikalisch-Mathematische Klasse},
  language  = {German},
  pages     = {144--153},
  pagetotal = {12 Seiten},
  ppn_gvk   = {1032800631},
  zbl       = {57.1147.01},
}

\end{document}